\documentclass[11pt,reqno]{amsart}
\usepackage{amssymb,amscd,amsbsy,mathrsfs}

\newcommand{\lb}{\linebreak}

\renewcommand{\a}{\alpha}
\renewcommand{\b}{\beta}

\renewcommand{\d}{\delta}
\newcommand{\e}{\varepsilon}

\newcommand{\z}{\zeta}

\renewcommand{\l}{\lambda}

\newcommand{\s}{\sigma}

\newcommand{\f}{\varphi}
\renewcommand{\o}{\omega}

\newcommand{\D}{\Delta}
\renewcommand{\L}{\Lambda}

\renewcommand{\O}{\Omega}

\newcommand{\B}{{\mathscr B}}

\newcommand{\cd}{{\mathscr D}}
\newcommand{\F}{{\mathscr F}}
\newcommand{\Q}{{\mathscr Q}}
\newcommand{\h}{{\mathscr H}}

\newcommand{\C}{{\Bbb C}}
\newcommand{\T}{{\Bbb T}}

\newcommand{\R}{{\Bbb R}}
\newcommand{\Z}{{\Bbb Z}}

\newcommand{\0}{{\boldsymbol{0}}}

\newcommand{\bs}{\boldsymbol}

\newcommand{\bS}{{\boldsymbol S}}

\newcommand{\rf}[1]{(\ref{#1})}

\newcommand{\df}{\stackrel{\mathrm{def}}{=}}

\newcommand{\supp}{\operatorname{supp}}

\newcommand{\rank}{\operatorname{rank}}
\newcommand{\const}{\operatorname{const}}

\newcommand{\eeq}{\end{equation}}
\newcommand{\beq}{\begin{equation}}
\newcommand{\bay}{\begin{eqnarray}}
\newcommand{\ba}{\begin{align*}}
\newcommand{\ea}{\end{align*}}
\newcommand{\ey}{\end{eqnarray}}
\newcommand{\bey}{\begin{eqnarray*}}
\newcommand{\eey}{\end{eqnarray*}}

\newcommand{\be}{\infty}

\newcommand{\bl}{\blacksquare}

\newcommand{\Pf}{{\bf Proof. }}
\newcommand{\im}{\operatorname{Im}}

\newtheorem{thm}{\hspace{\parindent}Theorem}[section]

\newtheorem{cor}[thm]{\hspace{\parindent}Corollary}
\newtheorem{lem}[thm]{\hspace{\parindent}Lemma}

\theoremstyle{remark}

\newtheorem*{rem*}{Remark}

\newcommand\sgn{\operatorname{sgn}}
\newcommand{\Bbbone}{{\rm{1\mathchoice{\kern-0.25em}{\kern-0.25em}{\kern-0.2em}{\kern-0.2em}I}}}
\newcommand{\fI}{{\frak I}}
\newcommand{\qm}{\quad\mbox{and}\quad}

\begin{document}

\newcommand{\vse}{\vspace{.2in}}
\numberwithin{equation}{section}

\title{Functions of perturbed unbounded self-adjoint operators. Operator Bernstein type inequalities}
\author{A.B. Aleksandrov and V.V. Peller}
\thanks{The first author is partially supported by RFBR grant 08-01-00358-a and by
Russian Federation presidential grant NSh-2409.2008.1;
the second author is partially supported by NSF grant DMS 0700995 and by ARC grant}

\newcommand{\mt}{{\mathcal T}}

\begin{abstract} 
This is a continuation of our papers \cite{AP2} and \cite{AP3}.
In those papers we obtained estimates for finite differences $(\D_Kf)(A)=f(A+K)-f(A)$ of the order $1$ and
$(\D_K^mf)(A)\df\sum\limits_{j=0}^m(-1)^{m-j}\left(\begin{matrix}m\\j\end{matrix}\right)f\big(A+jK\big)$
of the order $m$
for certain classes of functions $f$, where $A$ and $K$ are bounded self-adjoint operator.
In this paper we extend results of \cite{AP2} and \cite{AP3} to the case
of unbounded self-adjoint operators $A$.
Moreover, we obtain operator Bernstein type inequalities
for entire functions of exponential type.
This allows us to obtain alternative proofs of the main results of \cite{AP2}.
We also obtain operator Bernstein type inequalities for functions of unitary operators.
Some results of this paper as well as of the papers \cite{AP2} and \cite{AP3}
were announced in \cite{AP1}.
\end{abstract}

\maketitle

\

\begin{center}
{\Large Contents}
\end{center}

\

\begin{enumerate}
\item[1.] Introduction \quad\dotfill \pageref{s1}
\item[2.] Function spaces \quad\dotfill \pageref{s2}
\item[3.] Symmetrically normed and quasinormed ideals of operators\quad\dotfill \pageref{s3.5}
\item[4.] Moduli of continuity and unbounded self-adjoint operators  \quad\dotfill \pageref{s3}
\item[5.] Operator Bernstein type inequalities and their applications \quad\dotfill \pageref{s4}
\item[6.] H\"older--Zygmund estimates for self-adjoint operators \quad\dotfill \pageref{s5}
\item[7.] Higher order moduli of continuity \quad\dotfill \pageref{s6}
\item[] References \quad\dotfill \pageref{bibl}
\end{enumerate}

\

\setcounter{section}{0}
\section{\bf Introduction}
\setcounter{equation}{0}
\label{s1}

\

In this paper we are continuing our study of properties of functions of perturbed operators. It can be considered as a continuation of our papers \cite{AP2} and \cite{AP3}.

Recall that a Lipschitz function $f$ on the real line (i.e., a function $f$ satisfying the estimate
$$
|f(x)-f(y)|\le\const|x-y|,\quad x,~y\in\R,)
$$
does not have to satisfy the inequality
$$
\|f(A)-f(B)\|\le\const\|A-B\|
$$
for arbitrary self-adjoint operators $A$ and $B$ on Hilbert space, i.e.,
it does not have to be an {\it operator Lipschitz function}.
The existence of such functions was proved for the first time in \cite{F1}.
Later Kato established in \cite{Ka} that the function $t\mapsto|t|$ is not operator Lipschitz.

In \cite{Pe1} necessary conditions were found for a function $f$ to be operator Lipschitz.
Those necessary conditions also imply that Lipschitz  functions do not have to be operator Lipschitz.
In particular, it was shown in \cite{Pe1} that an operator Lipschitz function must belong locally to the Besov space $B_1^1(\R)$ (see \S\,2 for an introduction to Besov spaces). Note that in \cite{Pe1} and \cite{Pe3} a stronger necessary condition was also obtained.

On the other hand it was shown in \cite{Pe3} that if $f$ belongs to the Besov space $B_{\be1}^1(\R)$, then $f$ is operator Lipschitz.

For functions $f$ in the H\"older class $\L_\a(\R)$, $0<\a<1$, i.e., for functions satisfying the condition
$$
|f(x)-f(y)|\le\const|x-y|^\a,
$$
it was shown in \cite{F1} that
$$
\|f(A)-f(B)\|\le\const\|f\|_{\L_\a(\R)}\left(\log_2^2\frac{b-a}{\|A-B\|}+1\right)^\a\|A-B\|^\a,
$$
whenever $A$ and $B$ self-adjoint operators with spectra in $[a,b]$ (see also \cite{F2}).

For almost 40 years it was not known whether one can remove the logarithmic factor on the right-hand side of this inequality. In other words, it was unclear whether a H\"older function of order $\a$, $0<\a<1$, must be {\it operator
H\"older of order $\a$}, i.e.,
$$
\|f(A)-f(B)\|\le\const\|A-B\|^\a
$$
for self-adjoint operators $A$ and $B$ on Hilbert space.
Many mathematicians working on estimates of functions of perturbed operators believed that the answer should be negative.

It turned out, however, that the situation differs dramatically from the situation with Lipschitz functions. We proved in \cite{AP2} (see also \cite{AP1} where the main results of \cite{AP2} were announced) that H\"older functions of order $\a$, $0<\a<1$, {\it must be operator H\"older of order $\a$}. Note that Farforovskaya and Nikolskaya found in \cite{FN} a different proof of this result.

We obtained in \cite{AP2} sharp estimates for the norms of $f(A)-f(B)$ in terms of the norm of $A-B$
for various classes of functions $f$. Here $A$ and $B$ are self-adjoint operators on Hilbert space and
$f$ is a function on the real line $\R$.
We also obtained in \cite{AP2} sharp estimates for the norms of higher order differences
\bay
\label{hod}
\big(\D_K^mf\big)(A)\df
\sum_{j=0}^m(-1)^{m-j}\left(\begin{matrix}m\\j\end{matrix}\right)f\big(A+jK\big),
\ey
where $A$ and $K$ are self-adjoint operators.

In \cite{AP3} we found  sharp estimates for the Schatten--von Neumann norms of
first order differences $f(A)-f(B)$ and higher order differences
$\big(\D_K^mf\big)(A)$ for functions $f$ that belong to a H\"older--Zygmund class
$\L_\a(\R)$, $0<\a<\be$, (see \S\,\ref{s2} for the definition of these spaces).

In \cite{AP2} and \cite{AP3} we considered in detail the case
of arbitrary self-adjoint operators to estimate first order differences.
However, to estimate higher order differences, we gave proofs for bounded
self-adjoint operators.
In this paper we prove that the results of \cite{AP2} and \cite{AP3} are true also in the case of unbounded self-adjoint operators.

In \cite{AP2} we introduced four operator  moduli of continuity
$\O_f,\O_f^{[1]},\O_f^{[2]},\O_f^{[3]}$ for a function $f:\R\to\C$.
In \S\,\ref{s3} of this paper we continue studying operator moduli of continuity.
In particular, we show that in the definitions of the operator moduli of continuity we can allow unbounded self-adjoint operators.

In \S\,\ref{s4} we give sharp estimates of the finite differences
$\big(\D_K^mf\big)(A)$ for the class of all entire functions $f$
of exponential type $\s$ with $|f|\le1$ on $\R$. The proofs are elementary.
The previous proofs of similar results (see \cite{Pe1}, \cite{Pe3}, \cite{AP2}) are based on
techniques of multiple operator integrals. Besides, they do not allow us
to obtain estimates with best possible constants.
The results \S\,\ref{s4} allow us to obtain more elementary proofs of the main results of \cite{AP2}.

In \cite{AP2} we proved that $\O_f\le\O_f^{[1]}=\O_f^{[2]}=\O_f^{[3]}\le2\O_f$.
Our estimates of \S\,\ref{s4} allow us to prove that in general
$\O_f\not=\O_f^{[j]}$ for $j=1,2,3$.

In \S\,\ref{s5} we show that the results of \S\,4 of \cite{AP2}
and \S\,5, \S\,11 of  \cite{AP3}
are true also for unbounded self-adjoint operators, while
in \S\,\ref{s6} we extend the results of \S\,7 and \S\,11 of \cite{AP2}
to the case of unbounded self-adjoint operators.

In \S\,2 we collect necessary information on function classes, while in \S\,3 we give a brief introduction in operator ideals.

The authors are grateful to N.A. Shirokov for a question that
he asked during a seminar talk. Theorem \ref{opbern} below gives
a positive answer to this question.

\

\section{\bf Function spaces}
\setcounter{equation}{0}
\label{s2}

\

{\bf 2.1. Besov classes.}
The purpose of this subsection is to give a brief introduction to Besov spaces that
play an important role in problems of perturbation theory.
In this paper we consider Besov spaces on the real line.

Let $w$ be an infinitely differentiable function on $\R$ such
that
\bay
\label{w}
w\ge0,\quad\supp w\subset\left[\frac12,2\right],\quad\mbox{and}
\quad w(x)=1-w\left(\frac x2\right)\quad\mbox{for}\quad x\in[1,2].
\ey

Consider the functions $W_n$ and $W_n^\sharp$ in the Schwartz class $\mathscr S(\R)$ defined by
$$
\F W_n(x)=w\left(\frac{x}{2^n}\right),\quad\F W^\sharp_n(x)=\F W_n(-x),\quad n\in\Z,
$$
where $\F$ is the {\it Fourier transform}:
$$
\big(\F f\big)(t)=\int_\R f(x)e^{-{\rm i}xt}\,dx,\quad f\in L^1.
$$

With every tempered distribution $f\in{\mathscr S}^\prime(\R)$ we
associate a sequences $\{f_n\}_{n\in\Z}$,
$$
f_n\df f*W_n+f*W_n^\sharp.
$$
Initially we define the (homogeneous) Besov class $\dot B^s_{pq}(\R)$ as the set of all $f\in{\mathscr S}^\prime(\R)$
such that
\bay
\label{Wn}
\{2^{ns}\|f_n\|_{L^p}\}_{n\in\Z}\in\ell^q(\Z).
\ey
According to this definition, the space $\dot B^s_{pq}(\R)$ contains all polynomials.
Moreover, the distribution $f$ is defined by the sequence $\{f_n\}_{n\in\Z}$
uniquely up to a polynomial. It is easy to see that the series $\sum_{n\ge0}f_n$ converges
in ${\mathscr S}^\prime(\R)$.
However, the series $\sum_{n<0}f_n$ can diverge in general. It is easy to prove that the
series $\sum_{n<0}f_n^{(r)}$ converges on uniformly $\R$ for each nonnegative integer
$r>s-1/p$. Note that in the case $q=1$ the series $\sum_{n<0}f_n^{(r)}$ converges uniformly,
whenever $r\ge s-1/p$.

Now we can define the modified (homogeneous) Besov class $B^s_{pq}(\R)$. We say that a distribution $f$
belongs to $B^s_{pq}(\R)$ if $\{2^{ns}\|f_n\|_{L^p}\}_{n\in\Z}\in\ell^q(\Z)$ and
$f^{(r)}=\sum_{n\in\Z}f_n^{(r)}$ in the space ${\mathscr S}^\prime(\R)$, where $r$ is
the minimal nonnegative integer such that $r>s-1/p$ ($r\ge s-1/p$ if $q=1$). Now the
function $f$ is determined uniquely by the sequence $\{f_n\}_{n\in\Z}$ up
to a polynomial of degree less that $r$, and a polynomial $\f$ belongs to $B^s_{pq}(\R)$
if and only if $\deg\f<r$.

To define a regularized de la Vall\'ee Poussin type kernel $V_n$, we define the $C^\be$ function $v$ on $\R$ by
\bay
\label{VP}
v(x)=1\quad\mbox{for}\quad x\in[-1,1]\quad\mbox{and}\quad v(x)=w(|x|)\quad\mbox{if}\quad |x|\ge1,
\ey
where $w$ is a function described in \rf{w}.
We define de la Vall\'ee Poussin type functions $V_n$, $n\in\Z$, by
$$
\F V_n(x)=v\left(\frac{x}{2^n}\right),
$$
where $v$ is a function given by \rf{VP}.

In this paper an important role is played by
H\"older--Zygmund classes $\L_\a(\R)\df B^\a_{\be}(\R)$, $\a>0$. The class $\L_\a(\R)$ can be described
as the class of continuous functions $f$ on $\R$ such that
$$
\big|(\D^m_tf)(x)\big|\le\const|t|^\a,\quad t\in\R,
$$
where the difference operator $\D_t$ is defined by
$$
(\D_tf)(x)=f(x+t)-f(x),\quad x\in\R,
$$
and $m\in\Z$, $m-1\le\a<m$.

We can introduce the following equivalent (semi)norm on $\L_\a(\R)$:
$$
\sup_{n\in\Z}2^{n\a}\big(\|f*W_n\|_{L^\be}+\|f*W_n^\sharp\|_{L^\be}\big),\quad f\in\L_\a(\R).
$$

Consider now the class $\l_\a(\R)$, which is defined as the closure of the Schwartz class
${\mathscr S}(\R)$ in $\L_\a(\R)$.
The following result gives a description of $\l_\a(\R)$
for $\a>0$. We use the following notation: $C_0(\R)$ stands for the space of continuous functions $f$ on $\R$
such that $\lim\limits_{|x|\to\be}f(x)=0$.

\begin{thm}
\label{mH}
Let $\a>0$ and let $m$ be the integer such that $m-1\le\a<m$. Suppose that $f\in\L_\a(\R)$.
The following are equivalent:

{\em(i)} $f\in\l_\a(\R)$;

{\em(ii)} $f_n\in C_0(\R)$ for every $n\in\Z$ and
$$
\lim_{|n|\to\be}2^{n\a}\|f_n\|_{L^\be}=0;
$$

{\em(iii)} the following equalitites hold:
$$
\lim_{t\to0}|t|^{-\a}\big(\D_t^mf\big)(x)=0\quad\mbox{uniformly in}\quad x\in\R,
$$
$$
\lim_{|t|\to\be}|t|^{-\a}\big(\D_t^mf\big)(x)=0\quad\mbox{uniformly in}\quad x\in\R,
$$
and
$$
\lim_{|x|\to\be}|t|^{-\a}\big(\D_t^mf\big)(x)=0\quad\mbox{uniformly in}\quad t\in\R\setminus\{0\}.
$$
\end{thm}

We refer the reader to \cite{AP2} for the proof of this theorem.

The dual space $\big(\l_\a(\R)\big)^*$ to $\l_\a(\R)$ can be identified in a
natural way with $B_1^{-\a}(\R)$ with respect to the pairing
$$
\langle f,g\rangle\df \lim_{N\to\be}\sum_{n=-N}^N\int_\R\big( \F(f_n)\big)(t)\big(\F g\big)(t)\,dt,
\quad f\in\l_\a(\R),~g\in B_1^{-\a}(\R).
$$
The dual space $\big(B_1^{-\a}(\R)\big)^*$ to $B_1^{-\a}(\R)$ can be identified
with $\L_\a(\R)$ with respect to the same pairing.

We refer the reader to \cite{Pee} and \cite{Pe4} for more detailed information on Besov spaces.

\medskip

{\bf 2.2. Spaces $\bs{\L_\o}$.} Let $\o$ be a modulus of continuity, i.e., $\o$ is
a nondecreasing continuous function on $[0,\be)$
such that $\o(0)=0$, $\o(x)>0$ for $x>0$, and
$$
\o(x+y)\le\o(x)+\o(y),\quad x,~y\in[0,\be).
$$
We denote by $\L_\o(\R)$ the space of functions on $\R$ such that
$$
\|f\|_{\L_\o(\R)}\df\sup_{x\ne y}\frac{|f(x)-f(y)|}{\o(|x-y|)}.
$$

\begin{thm}
\label{Vn}
There exists a constant $c>0$ such that for an arbitrary
modulus of continuity $\o$ and for an arbitrary function $f$ in $\L_\o(\R)$,
the following inequality holds:
\bay
\label{VPn}
\|f-f*V_n\|_{L^\be}\le c\,\o\big(2^{-n}\big)\|f\|_{\L_\o(\R)},\quad n\in\Z.
\ey
\end{thm}

We refer the reader to \cite{AP2} for the proof of this theorem.

\begin{cor}
\label{Wnn}
Let $f\in\L_\o(\R)$. Then
$$
\|f*W_n\|_{L^\be}\le\const\o\big(2^{-n}\big)\|f\|_{\L_\o(\R)},\quad n\in\Z,
$$
and
$$
\|f*W^\sharp_n\|_{L^\be}\le\const\o\big(2^{-n}\big)\|f\|_{\L_\o(\R)},\quad n\in\Z.
$$
\end{cor}

\medskip

{\bf 2.3. Spaces $\bs{\L_{\o,m}}$.}
We proceed now to moduli of continuity of higher order. For a continuous function $f$ on $\R$,
we define the $m$th modulus of continuity $\o_{f,m}$ of $f$ by
$$
\o_{f,m}(x)=\sup_{\{h:0\le h\le x\}}\big\|\D_h^mf\big\|_{L^\be}=
\sup_{\{h:0\le|h|\le x\}}\big\|\D_h^mf\big\|_{L^\be},\quad x>0.
$$

The following elementary formula  can easily be verified by induction:
\bay
\label{ind}
\big(\D_{2h}^mf\big)(x)=\sum\limits_{j=0}^m{m\choose j}\big(\D_{h}^mf\big)(x+jh).
\ey
It follows from \rf{ind} that
$\o_{f,m}(2x)\le 2^m\o_{f,m}(x),\quad x>0$.

Suppose now that
$\o$ is a nondecreasing function on $(0,\be)$ such that
\bay
\label{on}
\lim_{x\to0}\o(x)=0\quad\mbox{and}\quad
\o(2x)\le2^m\o(x)\quad\mbox{for}\quad x>0.
\ey
 It is easy to see that in this case
 \bay
 \label{udv}
 \o(tx)\le2^mt^m\o(x),\quad\mbox{for all}\quad x>0\quad\mbox{and}\quad t>1.
 \ey
Denote by $\L_{\o,m}(\R)$ the set of continuous functions $f$ on $\R$
satisfying
$$
\|f\|_{\L_{\o,m}(\R)}\df\sup\limits_{t>0}\frac{\|\D^m_tf\|_{L^\infty}}{\o(t)}<+\infty.
$$

\begin{thm}
\label{mnn}
There exists $c>0$ such that for an arbitrary
nondecreasing function $\o$ on $(0,\be)$ satisfying {\em\rf{on}} and
for an arbitrary function $f\in\L_{\o,m}(\R)$, the following inequality holds:
$$
\|f-f*V_n\|_{L^\be}\le c\,\o\big(2^{-n}\big)\|f\|_{\L_{\o,m}(\R)},\quad n\in\Z.
$$
\end{thm}

We refer the reader to \cite{AP2} for the proof of this theorem.

\begin{cor}
\label{Wnm}
Let $f\in\L_{\o,m}(\R)$. Then
$$
\|f*W_n\|_{L^\be}\le\const\o\big(2^{-n}\big)\|f\|_{\L_\o(\R)},\quad n\in\Z,
$$
and
$$
\|f*W^\sharp_n\|_{L^\be}\le\const\o\big(2^{-n}\big)\|f\|_{\L_\o(\R)},\quad n\in\Z.
$$
\end{cor}

\

\

\section{\bf Symmetrically normed and quasinormed ideals of operators}

\setcounter{equation}{0}
\label{s3.5}

\

In this section we give a brief introduction to ideals of operators on Hilbert space.

 First, we remind the definition of singular values of bounded linear operators on Hilbert space. Let $T$ be a bounded linear operator. The singular values $s_j(T)$, $j\ge0$, are defined by
$$
s_j(T)=\inf\big\{\|T-R\|:~\rank R\le j\big\}.
$$
Clearly, $s_0(T)=\|T\|$, and $T$ is compact if and only if $s_j(T)\to0$ as $j\to\be$.

\medskip

{\bf Definition.} Let $\h$ be a Hilbert space and let ${\frak I}$ be a nonzero linear manifold in the set
$\B$=$\B(\h)$ of bounded linear operators on $\h$ that is equipped with a quasi-norm $\|\cdot\|_{\frak I}$ that makes
$\fI$ a quasi-Banach space.
We say that ${\frak I}$ is a {\it (symmetrically) quasinormed ideal} if for every $A$ and $B$ in $\B$ and
$T\in{\frak I}$,
\bay
\label{qni}
ATB\in{\frak I}\qm\|ATB\|_\fI\le\|A\|\cdot\|B\|\cdot\|T\|_\fI.
\ey
Note that if $\fI\ne\B$, then $\fI$ is contained in the set of all compact operators.
We put $\|T\|_\fI\df\be$ if $T\notin\fI$.
A quasinormed ideal $\fI$ is called a {\it normed ideal} if $\|\cdot\|_\fI$ is a norm.

Let $\bS_p$, $0<p<\be$, be the Schatten--von Neumann class of operators $T$ on Hilbert space
such that
$$
\|T\|_{\bS_p}\df\left(\sum_{j\ge0}\big(s_j(T)\big)^p\right)^{1/p}.
$$
This is a normed ideal for $p\ge1$. Denote by $\bS_\be$ the class of all compact operators on Hilbert space.
For $T\in\bS_\be$, we put $\|T\|_{\bS_\be}\df s_0(T)=\|T\|$.
Clearly, $\bS_1\subset\fI$ for every normed ideal $\fI$.

Let $l$ be a nonnegative integer and $p>0$. Put
$$
\|T\|_{\bS_p^l}\df\left(\sum_{j=0}^l\big(s_j(T)\big)^p\right)^{1/p}
$$
for $T\in\B$. Clearly, $\|T\|\le\|T\|_{\bS_p^l}\le(l+1)^{1/p}\|T\|$
for all operators $T$.
It is well known that $\|\cdot\|_{\bS_p^l}$ is a norm for $p\ge1$
(see \cite{BS}).

Note that if an operator $T$ is represented as the sum of two operators:
$T=T_1+T_2$, then $\|T\|_{\bS_1^l}\le\|T_1\|_{\bS_1^l}+\|T_2\|_{\bS_1^l}
\le\|T_1\|_{\bS_1}+(l+1)\|T_2\|$.
It is well known and it is easy to see that
\bay
\label{s1sbe}
\|T\|_{\bS_1^l}=\inf\Big\{\|T_1\|_{\bS_1}+(l+1)\|T_2\|:T=T_1+T_2\Big\},
\ey
see \cite{Mi}.
An analog of this formula for
symmetric spaces can be found in \cite{KPS}, Ch. 2, formula (3.5).

If $T_1,T_2\in\B$ and $s_j(T_2)\le s_j(T_1)$ for all $j\ge0$, then it follows from \rf{qni}
that the condition that $T_1\in\fI$ implies that $T_2\in\fI$ and $\|T_2\|_{\fI}\le\|T_1\|_{\fI}$ for every quasinormed ideal  $\fI$.
We say that a quasinormed ideal $\fI$ has {\it majorization property}  (respectively {\it weak majorization property}) if the conditions
$$
T_1\in\fI,\quad T_2\in\B,\quad
\mbox{and}\quad
\|T_2\|_{\bS_1^l}\le\|T_1\|_{\bS_1^l}\quad
\mbox{for all}\quad
l\ge0
$$
imply that
$$
T_2\in\fI\quad\mbox{and}\quad\|T_2\|_{\fI}\le\|T_1\|_{\fI}\quad(\text{respectively}\quad
\|T_2\|_{\fI}\le C\|T_1\|_{\fI})
$$
(see \cite{GK}).
Note that if a quasinormed ideal $\fI$ has weak majorization property, then we can introduce on it the following new equivalent quasinorm:
$$
\|T\|_{\widetilde\fI}\df\sup\{\|R\|_{\fI}:\|R\|_{\bS_1^l}\le\|T\|_{\bS_1^l}\,\,\,\text{for all}\,\,\,l\ge0\}
$$
such that
$(\fI,\|\cdot\|_{\widetilde\fI})$ has majorization property.

It is well known that every separable normed ideal and every normed ideal
that is dual to a separable normed ideal has majorization property, see \cite{GK}.
Clearly, $\bS_1\subset\fI$ for every quasinormed ideal $\fI$
with majorization property.
Note also that every quasinormed ideal $\fI$ with $\b_{\fI}<1$, where
$\b_{\fI}$ denote the upper Boyd index of $\fI$ (see, for example, \cite{AP3}
for the definition of the upper Boyd index),
has weak majorization property.

It is well known that every normed ideal with majorization property is
an interpolating Banach space between $\bS_1$ and $\B$.
The corresponding statement for symmetric space see in \cite{KPS},
Ch. 2, Theorem 4.2.

We need the following fact that generalizes the above
result on interpolation between $\bS_1$ and $\B$. Apparently it is known among experts.

\begin{thm}
\label{inter}
Let $\fI$ be a quasinormed ideal with majorization property
and let \lb$\frak A:\frak L\to\frak L$ be a linear transformation
on a linear subset $\frak L$ of $\B$ such that $\frak L\cap\bS_1$
is dense in $\bS_1$. Suppose that
$\|\frak A T\|\le\|T\|$ and $\|\frak A T\|_{\bS_1}\le\|T\|_{\bS_1}$
for all $T\in\frak L$. Then $\|\frak A T\|_{\fI}\le\|T\|_{\fI}$
for every  $T\in\frak L$.
\end{thm}
\Pf The identity \rf{s1sbe} implies that $\|\frak A T\|_{\bS_1^l}\le\|T\|_{\bS_1^l}$
for all $T\in\frak L$ and $l\ge0$. Hence, $\|\frak A T\|_{\fI}\le\|T\|_{\fI}$
for every quasinormed ideal with majorization property and every $T\in\frak L$.
$\bl$

\begin{cor}
Under the hypotheses of Theorem {\em\ref{inter}},
$$
\|\frak A T\|_{\bS_p}\le\|T\|_{\bS_p}
$$
for all $p\ge1$ and $T\in\frak L$.
\end{cor}

The {\it Schur product} of matrices $C=\{c_{jk}\}_{j,k\ge0}$ and
$D=\{d_{jk}\}_{j,k\ge0}$ is defined as the matrix $C\star D$ whose entries
are the products of the entries of $C$ and $D$:
$$
C\star D=\{c_{jk}d_{jk}\}_{j,k\ge0}.
$$
Here we identify bounded linear operators on $\ell^2$
with their matrix representations with respect to the standard orthonormal basis of
$\ell^2$.
We denote by $\|C\|$ and $\|C\|_{\bS_1}$ the operator norm and the trace norm of the corresponding
operator on $\ell^2$ and we say that $\|C\|\df\be$ (respectively $\|C\|_{\bS_1}\df\be$)  if the matrix $u$
does not determine a bounded operator (respectively an operator of class $\bS_1$).
Finally, we use {\it the notation ${\frak C}^{00}(\Z_+^2)$ for the class of matrices $C=\{c_{jk}\}$ such that
the set $\{(j,k):c_{jk}\not=0\}$ is finite.}

We need the following known result that follows from the fact that the dual spaces
$(\bS_\be)^*$ and $(\bS_1)^*$ can naturally be identified with
$\bS_1$ and $\B$.

\begin{thm}
\label{mult}
Let $M$ be a matrix. Then
\begin{align*}
\|M\|_{\frak M}&\df\sup\big\{\|M\star C\|:~\|C\|=1\big\}\\[.2cm]
&=
\sup\big\{\|M\star C\|:~C\in{\frak C}^{00}(\Z_+^2),~\|C\|=1\big\}\\[.2cm]
&=\sup\big\{\|M\star C\|_{\bS_1}:~\|C\|_{\bS_1}=1\big\}\\[.2cm]
&=\sup\big\{\|M\star C\|_{\bS_1}:~C\in{\frak C}^{00}(\Z_+^2),~\|C\|_{\bS_1}=1\big\}.
\end{align*}
\end{thm}

A matrix $M$ is said to be a {\it Schur multiplier} if $\|M\|_{\frak M}<\be$.
 Denote by $\frak M$ the set of all Schur multipliers.

The following theorem is well known.

\begin{thm}
Let $\fI$ be a quasinormed ideal with majorization property.
Then
$$
\|M\star C\|_{\fI}\le\|M\|_{\frak M}\cdot\|C\|_{\fI}
$$
for every matrices $M$ and $C$. In particular,
$$
\|M\star C\|_{\bS_p}\le\|M\|_{\frak M}\cdot\|C\|_{\bS_p},\quad p\ge1,
$$
for all matrices $M$ and $C$.
\end{thm}

\Pf The result readily follows from Theorem \ref{inter}. $\bl$

%

%

%

\

\section{\bf Moduli of continuity and unbounded self-adjoint operators}
\setcounter{equation}{0}
\label{s3}

\

In this section we study properties of various operator moduli of continuity.
In particular, we show in this section that in the definition of operator moduli of continuity given in \cite{AP2} one can allow unbounded self-adjoint operators.

Let $f$ be a continuous function on $\R$.
We considered in \cite{AP2} the following four versions of operator moduli of continuity of $f$
that are defined on $(0,\be)$:
\begin{align*}
\O_f(\d)&\df\sup\big\{\|f(A)-f(B)\|:~A=A^*,~B=B^*,~\|A-B\|<\d\big\};\\[.2cm]
\O_f^{[1]}(\d)&\df\sup\big\{\|f(A)R-Rf(A)\|:
~A=A^*,~R=R^*,~\|R\|=1,~\|AR-RA\|<\d\big\};\\[.2cm]
\O_f^{[2]}(\d)&\df\sup\big\{\|f(A)R-Rf(A)\|:~A=A^*,~\|R\|=1,~\|AR-RA\|<\d\big\};\\[.2cm]
\O_f^{[3]}(\d)&\df\sup\big\{\|f(A)R-Rf(B)\|:~A=A^*,~B=B^*,
~\|R\|=1,~\|AR-RB\|<\d\big\}.
\end{align*}

In these definitions we assume that the operators $A$ and $B$ are bounded.

The following inequalities hold:
$$
\O_f\le\O_f^{[1]}=\O_f^{[2]}=\O_f^{[3]}\le2\O_f,
$$
see \cite{AP2}, Theorem 10.2.

Put $\O_f^{\flat}\df\O_f^{[1]}=\O_f^{[2]}=\O_f^{[3]}$.
We show in \S\,\ref{s4} that in general $\O_f\not=\O_f^{\flat}$.

\begin{thm}
\label{41}
Let $f$ be a continuous function on $\R$. Then the function
\bay
\label{do}
\d\mapsto\d^{-1}\O_f^\flat(\d),\quad\d>0,
\ey
is nonincreasing. In particular,
$$
\O_f^\flat(\d_1+\d_2)\le\O_f^\flat(\d_1)+\O_f^\flat(\d_2),\quad\d_1,~\d_2>0.
$$
\end{thm}

\Pf It suffices to verify that $\O_{f}^{\flat}(\d/\tau)\le\tau^{-1}\O_{f}^{\flat}(\d)$
for $\d\in(0,\be)$ and $\tau\in(0,1)$.
We have
\begin{align*}
\O_{f}^{\flat}(\d/\tau)&=\sup\big\{\|f(A)R-Rf(A)\|:A=A^*, ~\|R\|=1,~\|AR-RA\|<\d/\tau\big\}\\[.2cm]
&=\tau^{-1}\sup\big\{\|f(A)R-Rf(A)\|:A=A^*,~\|R\|=\tau,~\|AR-RA\|<\d\big\}.
\quad
\end{align*}
Obviously, for every operator $R$ with $\|R\|\le1$, there exists $\l\in\R$
such that $\|R+\l I\|=1$. Taking into account the fact that
$X(R+\l I)-(R+\l I)X=XR-RX$ for every operator $X$, we obtain
\begin{align}
\label{a-1}
&\sup\big\{\|f(A)R-Rf(A)\|:A=A^*,~\|R\|=\tau,~\|AR-RA\|<\d\big\}\nonumber\\
\le&\sup\big\{\|f(A)R-Rf(A)\|:A=A^*,~\|R\|=1,~\|AR-RA\|<\d\big\}=\O_{f}^{\flat}(\d).
\end{align}
Now the desired inequality is evident. $\bl$

\begin{cor}
The function
$\O_f^\flat$ is continuous as a function from $(0,\be)$ to $[0,\be]$.
\end{cor}

\Pf It suffices to observe that the function $\O_f^\flat(\delta)$ is nondecreasing and the function
\rf{do} is nonincreasing. $\bl$

\begin{cor}
\label{123cor}
In the definition of $\O_{f}^{[j]}$, $j=1,2,3$,
one can replace the condition $\|R\|=1$ with the condition $\|R\|\le1$.
\end{cor}

\Pf The case $j=2$ follows from inequality \rf{a-1}. A similar argument also
works for $j=1$. Let $j=3$. Then for $\tau\in(0,1)$, we have
\begin{align*}
\sup\big\{\|f(A)R&-Rf(B)\|:A=A^*,~B=B^*,~\|R\|=\tau,~\|AR-RB\|<\d\big\}\\
&=\tau\Omega_f^{[3]}(\d/\tau)=\tau\Omega_f^{\flat}(\d/\tau)
\le\Omega_f^{\flat}(\d)=\Omega_f^{[3]}(\d).\quad\bl
\end{align*}

%
%

\medskip

{\bf Remark.} It is easy to see that $\O_f(\d_1+\d_2)\le\O_f(\d_1)+\O_f(\d_2)$.
Hence, $\O_f$ is continuous if $\lim_{\d\to0}\O_f(\d)=0$.
However, we do not know whether the function
$$
\d\mapsto\d^{-1}\O_f(\d),\quad\d>0,
$$
is nonincreasing.

\medskip

It was shown in \cite{AP2},
Th.\,\,8.3, that if we allow unbounded
self-adjoint operators $A$ and $B$ in the definition of the operator modulus of continuity $\O_f$,
we obtain the same operator modulus of continuity $\O_f$.
In this section we prove that the same is true for $\O_f^{[1]}$, $\O_f^{[2]}$
and $\O_f^{[3]}$. Let us explain what we mean by $\|f(A)R-Rf(B)\|$
for not necessarily bounded self-adjoint operators $A$ and $B$.
Note that the operators $f(A)$ and $f(B)$ are normal.

Let $M$ and $N$ be (not necessarily bounded) normal operators
in a Hilbert space and let $R$ be a bounded operator on the same Hilbert space.
We say that the {\it operator $MR-RN$ is bounded} if
$R(\cd_N)\subset \cd_M$ and $\|MRu-RNu\|\le C\|u\|$ for every $u\in \cd_N$.
Then there exists a unique bounded operator $K$ such that
$Ku=MRu-RNu$ for all $u\in \cd_N$. In this case we write $K=MR-RN$.
Thus $MR-RN$ is bounded if and only if
\bay
\label{MN}
\big|(Ru,M^*v)-(Nu,R^*v)\big|\le C\|u\|\cdot\|v\|
\ey
for every $u\in \cd_N$ and $v\in \cd_{M^*}=\cd_{M}$. It is easy to see that
$MR-RN$ is bounded if and only if $N^*R^*-R^*M^*$ is bounded,
and $(MR-RN)^*=-(N^*R^*-R^*M^*)$.
In particular, we write $MR=RN$ if $R(\cd_N)\subset \cd_M$ and $MRu=RNu$ for every $u\in \cd_N$.
We say that $\|MR-RN\|=\be$ if $MR-RN$ is not a bounded operator.

We need the following obvious observation.

\medskip

{\bf Remark.} Let $M$ and $N$ be normal operators. Suppose that $M^*$ is the closure of an operator
$M_\flat$ and $N$ is the closure of an operator $N_\sharp$. Suppose that inequality
\rf{MN} holds for all $u\in \cd_{N_\sharp}$ and $v\in\cd_{M_\flat}$.
Then it holds for all $u\in \cd_N$ and $v\in \cd_{M}$.

\begin{lem}
\label{anbnrn}
Let $A$ and $B$ be self-adjoint operators and let
$R$ be an operator of norm $1$. Then there exist
a sequence of operators $\{R_n\}_{n\ge1}$ and sequences of
bounded self-adjoint operators $\{A_n\}_{n\ge1}$ and $\{B_n\}_{n\ge1}$ such that

{\em(i)} the sequence $\{\|R_n\|\}_{n\ge1}$ is nondecreasing and
$$
\lim_{n\to\be}\|R_n\|=1;
$$

{\em(ii)}
$$
\lim_{n\to\be}R_n=R
$$
in the strong operator topology;

{\em(iii)} for every
continuous functions $f$ on $\R$, the sequence
$\big\{\big\|f(A_n)R_n-R_nf(B_n)\big\|\big\}_{n\ge1}$ is nondecreasing
and
$$
\lim_{n\to\be}\big\|f(A_n)R_n-R_nf(B_n)\big\|=\|f(A)R-Rf(B)\|;
$$

{\em(iv)} if $f$ is a continuous function on $\R$ such that
$\|f(A)R-Rf(B)\|<\be$, then
$$
\lim_{n\to\be}f(A_n)R_n-R_nf(B_n)=f(A)R-Rf(B)
$$
in the strong operator topology;

{\em(v)} if $f$ is a continuous function on $\R$ such that
$\|f(A)R-Rf(B)\|<\be$, then for every $j\ge0$, the sequence
$\big\{s_j\big(f(A_n)R_n-R_nf(B_n)\big)\big\}_{n\ge1}$ is nondecreasing and
$$
\lim_{n\to\be}s_j\big(f(A_n)R_n-R_nf(B_n)\big)=s_j\big(f(A)R-Rf(B)\big).
$$
\end{lem}

\Pf
Put $P_n\df E_A\big([-n,n])$ and $Q_n\df E_B\big([-n,n])$,
where $E_A$ and $E_B$ are the spectral measures of $A$ and $B$. Put $A_n\df P_nA=AP_n$ and
$B_n\df Q_nB=BQ_n$.
Clearly,
\bay
\label{pnqn}
P_n\big(f(A)R-Rf(B)\big)Q_n=f(A_n)P_nRQ_n-P_nRQ_nf(B_n),\quad n\ge1.
\ey
It remains to put $R_n\df P_nRQ_n$. $\bl$

\begin{lem}
\label{commut}
Suppose that $AR=RB$ for a bounded operator $R$ and self-adjoint operators $A$ and $B$.
Then $f(A)R=Rf(B)$ for every continuous  function $f$ on $\R$.
\end{lem}

\Pf This is well known if $A=B$. The general case reduces to this special case by considering the operators
$$
\left(\begin{matrix}A&0\\0&B\end{matrix}\right)\quad\mbox{and}\quad
\left(\begin{matrix}0&R\\0&0\end{matrix}\right).\quad\bl
$$

\begin{thm}
Let $A$ and $B$ be self-adjoint operators and let
$R$ be a bounded operator such that $\|R\|=1$. Suppose that $AR-RB$
is bounded.
Then for every continuous function $f$ on $\R$, the following inequality holds:
$$
\big\|f(A)R-Rf(B)\big\|\le\O_f^{\flat}\big(\|AR-RB\|\big),
$$
where $\O_f^{\flat}(0)\df0$.
\end{thm}

\Pf
Lemma \ref{commut} allows us to restrict ourselves to the case
when $\|AR-RB\|>0$.
We have
\begin{align*}
\|f(A)R-Rf(B)\|&=\lim_{n\to\be}\|f(A_n)R_n-R_nf(B_n)\|\\[.2cm]
&\le\lim_{n\to\be}\O_f^{\flat}\big(\|A_nR_n-R_nB_n\|\big)
=\O_f^{\flat}\big(\|AR-RB\|\big),
\end{align*}
where $A_n$, $B_n$ and $R_n$ are as in Lemma \ref{anbnrn}. $\bl$

\begin{cor}
If we allow unbounded self-adjoint operators $A$ and $B$ in the definitions of the operator moduli of continuity $\O_f^{[j]}$, $j=1,\,2,\,3$, the result will be the same.
\end{cor}

The following theorem was proved in \cite{AP2} (Th. 10.1) in the case of bounded operators $A$ and $B$,
see also \cite{KS}.

\begin{thm}
\label{usr}
Let $f$ be a continuous function on $\R$. The following are equivalent:

{\em(i)} $\|f(A)-f(B)\|\le\|A-B\|$ for arbitrary self-adjoint operators $A$ and $B$;

{\em(ii)} $\|f(A)-f(B)\|\le\|A-B\|$ for all pairs of unitarily equivalent self-adjoint operators
$A$ and $B$;

{\em(iii)} $\|f(A)R-Rf(A)\|\le\|AR-RA\|$ for arbitrary self-adjoint operators $A$ and $R$
with $\|R\|<\be$;

{\em(iv)} $\|f(A)R-Rf(A)\|\le\|AR-RA\|$ for all self-adjoint operators $A$ and bounded operators $R$;

{\em(v)} $\|f(A)R-Rf(B)\|\le\|AR-RB\|$ for arbitrary self-adjoint operators $A$ and $B$
and an arbitrary bounded operator $R$.
\end{thm}

The same reasoning show that this theorem remains valid for arbitrary (not necessary bounded)
self-adjoint operators $A$ and $B$.

\begin{cor} Let $f$ be a continuous function on $\R$. Then
$$
\sup_{t>0}\frac{\O_f(t)}t=\sup_{t>0}\frac{\O_f^{\flat}(t)}t.
$$
\end{cor}

It is easy to see that a function $f\in C(\R)$ is  operator Lipschitz
(see the Introduction) if and only if
$$
\sup_{t>0}\frac{\O_f(t)}t=\sup_{t>0}\frac{\O_f^{\flat}(t)}t<\be.
$$

Recall that Kato \cite{Ka} proved that the function $t\mapsto|t|$ is not operator Lipschitz.
In \cite{AP2} we noted that this result by Kato implies that $\O_{|x|}(\d)=\O_{|t|}^{\flat}(\d)=\be$ for $\d>0$,
see the example following Theorem 8.2 in \cite{AP2}. This implies the following result:

\begin{thm}
\label{dif+-}
Let $f$ be a continuous function on $\R$ such that
$$
\lim_{t\to+\be}\frac{f(t)}{t}\not=\lim_{t\to-\be}\frac{f(t)}{t}
$$
and both limits exist and are finite.
Then $\O_f(\d)=\be$ for all $\d>0$.
\end{thm}

\Pf It suffices to consider the case when
$$
\lim_{+\be}t^{-1}f(t)=1\quad\mbox{and}\quad
\lim_{-\be}t^{-1}f(t)=-1.
$$
Put $f_n(t)=n^{-1}f(nt)$.
Clearly, $\lim_{n\to\be} f_n(t)=|t|$ and by Theorem \ref{41},
$\O_{f_n}^{\flat}(\d)=n^{-1}\O_{f}^{\flat}(n\d)\le\O_{f}^{\flat}(\d)$
for every $n\ge1$. Hence,
$\O_{|x|}^{\flat}\le\O_{f}^{\flat}$. It remains to observe that
$\O_{|x|}^{\flat}=\be$. $\bl$

We obtain a stronger result in Theorem \ref{difbe}.

The following theorem follows essentially (at least up to a multiplicative
constant) from results of \cite{KS}. We give a proof here for the reader's convenience.

\begin{thm}
\label{fI}
Let $\fI$ be a quasinormed ideal with majorization property and
let $f$ be a function satisfying the equivalent statements of Theorem {\em\ref{usr}}.
Then
$$
\|f(A)R-Rf(B)\|_{\fI}\le\|AR-RB\|_{\fI}
$$
for arbitrary self-adjoint operators $A$ and $B$ and
an arbitrary bounded operator $R$.
\end{thm}

\begin{cor}
\label{lipsp}
Let $f$ be a function satisfying the equivalent statements of Theorem {\em\ref{usr}}.
Then
$$
\|f(A)R-Rf(B)\|_{\bS_p}\le\|AR-RB\|_{\bS_p}\,,\quad 1\le p\le\be,
$$
for arbitrary self-adjoint operators $A$ and $B$ and
an arbitrary bounded operator $R$.
\end{cor}

{\bf Proof of Theorem \ref{fI}.} First we assume that $A$ and $B$ are self-adjoint operators
with pure point spectra such that $\s_{\rm point}(A)\cap\s_{\rm point}(B)=\varnothing$.
Then there exist orthonormal bases $\{e_j\}_{j\ge0}$ and $\{e_k^{\,\prime}\}_{k\ge0}$
of eigenvectors of $A$ and $B$. Let $Ae_j=\l_j e_j$ and $Be_k^{\,\prime}=\mu_k e_k^{\,\prime}$.
We identify each bounded operator $T$
with the matrix $\big\{(Te_k^{\,\prime}, e_j)\big\}_{j,k\ge0}$.
Thus the operator $AR-RB$ is identified with
the matrix
\bay
\label{rmatr}
C_R=\{c_{jk}\}_{j,k\ge0}\df\big\{(AR-RB)e_k^{\,\prime},e_j)\big\}_{j,k\ge0}
=\big\{(\l_j-\mu_k)(Re_k^{\,\prime},e_j)\big\}_{j,k\ge0}.
\ey
Denote by $\frak L$ the linear span of the rank one operators of the form $(\,\cdot\,,e_k^{\,\prime})e_j$
with $j,k\ge0$. Clearly, that transformer $R\mapsto C_R$ maps $\frak L$ onto ${\frak C}^{00}(\Z_+^2)$.
The operator $f(A)R-Rf(B)$ has matrix
\begin{align*}
C_R^{(f)}&=\left\{c_{jk}^{(f)}\right\}_{j,k\ge0}
\df\big\{(f(A)R-Rf(B))e_k^{\,\prime},e_j)\big\}_{j,k\ge0}\\[.2cm]
&=\big\{(f(\l_j)-f(\mu_k))(Re_k^{\,\prime},e_j)\big\}_{j,k\ge0}=M\star C_R,
\end{align*}
where
$$
M=M(f)=\left\{\frac{f(\l_j)-f(\mu_k)}{\l_j-\mu_k}\right\}_{j,k\ge0}.
$$
The condition $\|f(A)R-f(B)R\|\le\|AR-RB\|$ for all bounded operators $R$ and Theorem \ref{mult}
imply that $\|M\|_{\frak M}\le1$. Hence, by Theorem \ref{inter}, we have
$\|M\star C\|_{\fI}\le\|C\|_{\fI}$ for all matrices $C=C_R$ of the form \rf{rmatr}
with $R\in\B$ and all quasinormed ideals $\fI$
with majorization property.
Thus
\bay
\label{fARfB}
\|f(A)R-Rf(B)\|_{\fI}\le\|AR-RB\|_{\fI}
\ey
for all $R\in\B$ and all quasinormed ideal $\fI$
with majorization property. In particular,
\bay
\label{s1l}
\|f(A)R-Rf(B)\|_{\bS_1^l}\le\|AR-RB\|_{\bS_1^l}
\ey
for all $R\in\B$ and all $l\ge0$.

Let now $A$ and $B$ be arbitrary self-adjoint operators. Then we can
construct two sequences of self-adjoint operators with pure point spectra
$\{A_n\}$ and $\{B_n\}$ such that $\s_{\rm point}(A_n)\cap\s_{\rm point}(B)=\varnothing$
for all $n$, $\|A_n-A\|\to 0$ and $\|B_n-B\|\to0$ as $n\to\be$.
We have
$$
\|f(A_n)R-Rf(B_n)\|_{\bS_1^l}\le\|A_nR-RB_n\|_{\bS_1^l}
$$
for all $n$ and $l\ge0$. Passing to the limit as $n\to\be$,
we obtain inequality \rf{s1l} for all self-adjoint operators $A$ and $B$.
It remains to observe that inequalities \rf{s1l} for $l\ge0$
imply inequality \rf{fARfB}, because $\fI$ has majorization property.
$\bl$

Now we state analogues of Theorems \ref{usr} and \ref{fI}
for the unitary operators.

\begin{thm}
\label{usru}
Let $f$ be a continuous function on the unit circle $\T$. The following are equivalent:

{\em(i)} $\|f(U)-f(V)\|\le\|U-V\|$ for arbitrary unitary operators $U$ and $V$;

{\em(ii)} $\|f(U)-f(V)\|\le\|U-V\|$ for all pairs of unitarily equivalent unitary operators
$U$ and $V$;

{\em(iii)} $\|f(U)R-Rf(U)\|\le\|UR-RU\|$ for every unitary operator $U$ and
a bounded self-adjoint operator $R$;

{\em(iv)} $\|f(U)R-Rf(U)\|\le\|UR-RU\|$ for every unitary operator $U$ and
a bounded operator $R$;

{\em(v)} $\|f(U)R-Rf(V)\|\le\|UR-RV\|$ for arbitrary unitary operators $U$ and $V$
and an arbitrary bounded operator $R$.
\end{thm}

\begin{thm}
\label{fIu}
Let $\fI$ be a quasinormed ideal with majorization property and
let $f$ be a function satisfying the equivalent statements of Theorem {\em\ref{usru}}.
Then
$$
\|f(U)R-Rf(V)\|_{\fI}\le\|UR-RV\|_{\fI}
$$
for arbitrary unitary operators $U$ and $V$ and
an arbitrary bounded operator $R$.
\end{thm}

\begin{cor}
\label{lipspu}
Let $f$ be a function satisfying the equivalent statements of Theorem {\em\ref{usru}}.
Then
$$
\|f(U)R-Rf(V)\|_{\bS_p}\le\|UR-RV\|_{\bS_p}\,,\quad 1\le p\le\be,
$$
for arbitrary unitary operators $U$ and $V$ and
an arbitrary bounded operator $R$.
\end{cor}

We omit the proofs of Theorems \ref{usru} and \ref{fIu} because they repeat word-by-word
the proofs of Theorems \ref{usr} and \ref{fI}.



Theorem 4.1 in \cite{JW} implies that
every operator Lipschitz function $f$
is differentiable at every point. It is well known that the same argument gives the differentiability at $\infty$
in the following sense: there exists a finite limit $\lim\limits_{|x|\to+\infty}x^{-1}f(x)$.

\begin{thm}
\label{difbe}
Let $f\in C(\R)$.
Suppose that $\O_f(\d)<\be$ for $\d>0$.
Then the limit
$$
\lim_{|t|\to\be}t^{-1}f(t)
$$
exists and is finite.
\end{thm}

\Pf Let $\o_f$ denote the usual (scalar) modulus of continuity of $f$. Note that $\o_f(\d)\le\O_f(\d)\le\O_f^\flat(\d)\le\d\O_f^\flat(1)$ for $\d\ge1$.
Hence, $\limsup\limits_{|t|\to \be}|t|^{-1}|f(t)|<\be$.
Assume that the limit $\lim_{|t|\to\be}t^{-1}f(t)$ does not exist.
Then, as we have observed, $f$ cannot be operator Lipschitz.

For reader's convenience we repeat the corresponding arguments of \cite{JW} (see also \cite{Mc})
which allows us to prove (together with the fact that $f$ is not an operator Lipschitz function)
that
$\O_f(\d)=\be$ for $\d>0$.
The function $t\mapsto t^{-1}f(t)$ has at least two limit points as $|t|\to\be$.
Without loss of generality we may assume that $1$ and $-1$ are two such
limit points.
Then there exist two sequences $\{\l_j\}_{j\ge0}$ and $\{\mu_j\}_{j\ge0}$ in $\R$ such that

a) $0<2^{j+1}|\l_j|<|\mu_j|$ and $2^{j+2}|\mu_j|<|\l_{j+1}|$ for all $j\ge0$;

b) $|f(\l_j)-\l_j|<2^{-j-1}|\l_j|$ and $|f(\mu_j)+\mu_j|<2^{-j-1}|\mu_j|$
for all $j\ge0$.

Put $m_{jk}\df\dfrac{f(\l_j)-f(\mu_k)}{\l_j-\mu_k}$.
Let $0\le j\le k$. Then
\bey
|m_{jk}+1|=\left|\frac{f(\l_j)+\l_j-f(\mu_k)-\mu_k}{\l_j-\mu_k}\right|
\le\frac{3|\l_k|+|f(\mu_k)+\mu_k|}{|\mu_k|-|\l_k|}\le3\cdot2^{-k}+2^{-k}=4\cdot2^{-k}.
\eey
If $0\le k<j$, then
\bey
|m_{jk}-1|=\left|\frac{f(\l_j)-\l_j-f(\mu_k)+\mu_k}{\l_j-\mu_{k}}\right|
\le\frac{|f(\l_j)-\l_j|+3|\mu_{j-1}|}{|\l_j|-|\mu_{j-1}|}\le3\cdot2^{-j}+2^{-j}=4\cdot2^{-j}.
\eey
Hence,
$$
\left|m_{jk}-\sgn\Big(j-k-\frac12\Big)\right|\le4\cdot2^{-\max(j,k)}\le4\cdot 2^{-\frac{j+k}2},
$$
whence $\big\{m_{jk}-\sgn(j-k-\frac12)\big\}_{j,k\ge0}\in\frak M$. It is well known
that $\{\sgn(j-k-\frac12)\}_{j,k\ge0}\not\in\frak M$. Thus
$\{m_{jk}\}_{j,k\ge0}\not\in\frak M$.

Let us consider diagonal self-adjoint operators $A$ and $B$ such that
$Ae_j=\l_je_j$ and $Be_j=\mu_j e_j$, where $\{e_j\}_{j\ge0}$ is an orthonormal
basis.
Let $R$ be a bounded operator such that $\{(Re_k,e_j)\}_{j,k\ge0}\in {\frak C}^{00}(\Z_+^2)$.
Then the operators $AR-RB$ and $f(A)R-Rf(B)$ are well defined, $((AR-RB)e_k,e_j)=(\l_j-\mu_k)(Re_k,e_j)$ and
$((f(A)R-Rf(B))e_k,e_j)=(f(\l_j)-f(\mu_k))(Re_k,e_j)$.

Since $\{m_{jk}\}_{j,k\ge0}\not\in\frak M$, it follows that
for every $M>0$, there exists an operator $R$ such that
$\|AR-RB\|=C_0$ and $\|f(A)R-Rf(B)\|>MC_0$, where $C_0$ is a positive number  that will be chosen later.

We have
$$
|\l_j-\mu_k|^{-1}\le\left\{\begin{array}{ll}2|\l_j|^{-1},&\text {if}\,\,\,\,j>k,\\[.2cm]
2|\mu_k|^{-1},&\text {if}\,\,\,\,j\le k.
\end{array}\right.
$$
Hence, $|\l_j-\mu_k|^{-1}\le C\cdot2^{-j-k}$ because $|\l_{j+1}|\cdot|\l_j|^{-1}\to\be$
and $|\mu_{j+1}|\cdot|\mu_j|^{-1}\to\be$ as \lb $j\to\be$.
Thus $\{(\l_j-\mu_k)^{-1}\}_{j,k\ge0}\in\frak M$.
Now it is clear that \lb $\|R\|\le C_0\big\|\{(\l_j-\mu_k)^{-1}\}_{j,k\ge0}\big\|_{\frak M}=1$
if we put $C_0\df\big\|\{(\l_j-\mu_k)^{-1}\}_{j,k\ge0}\big\|_{\frak M}^{-1}$.
Hence, $\O_f^\flat(C_0)=\be$, and we get a contradiction. $\bl$

\medskip

{\bf Remark.}
Let $f$ be a continuous function
defined on a closed subset $E$ of $\R$. We can define in a similar way the operator moduli of continuity
$\O_{f,E}$, $\O_{f,E}^{[1]}$, $\O_{f,E}^{[2]}$ and $\O_{f,E}^{[3]}$
of $f$ if we consider only self-adjoint operators $A$ and $B$
with $\s(A),\s(B)\subset E$. In the same way we can prove that
$$
\O_{f,E}\le\O_{f,E}^{[1]}=\O_{f,E}^{[2]}=\O_{f,E}^{[3]}\le2\O_{f,E}
$$
and the function $\d\mapsto\d^{-1}\O_{f,E}^{\flat}(\d)$ is nonincreasing,
where $\O_{f,E}^{\flat}\df\O_{f,E}^{[1]}=\O_{f,E}^{[2]}=\O_{f,E}^{[3]}$.

Almost all results on operator moduli of continuity can be
extend to this case. As before, to obtain the corresponding results for unbounded self-adjoint operators, we use the construction in the proof of Lemma \ref{anbnrn}. Let us observe that the operators $A_n$ and $B_n$ constructed there
satisfy the following conditions: $\s(A_n)\subset\{0\}\cup\s(A)$ and
$\s(B_n)\subset\{0\}\cup\s(B)$. Thus everything works in the same way in the case
$0\in E$. The general case reduces to this special case with the help of translations.


\medskip

Theorem \ref{usr} also admits a natural generalization to the case of functions $f$ defined on $E$, see also \cite{KS}.

Let us state the corresponding generalization of Theorem \ref{difbe}.
\begin{thm}
\label{difbe+}
Let $f$ be a continuous function on an unbounded closed subset $E$ of $\R$.
Suppose that $\O_{f, E}(\d)<\be$ for $\d>0$.
Then
the function $t\mapsto t^{-1}f(t)$ has a finite limit as $|t|\to\be$, $t\in E$.
\end{thm}




\

\section{\bf Operator Bernstein type inequalities and their applications}
\setcounter{equation}{0}
\label{s4}

\

The results of \cite{Pe1} imply that for a trigonometric polynomial $f$ of degree $d$ and unitary operators $U$ and $V$, the following inequality holds:
$$
\|f(U)-f(V)\|\le\const d\|f\|_{L^\be}\|U-V\|.
$$
On the other hand, the results of \cite{Pe3} imply that if $f$ is and entire function of exponential type at most $\s$ that is bounded on $\R$, then for arbitrary self-adjoint operators $A$ and $B$, the following inequality holds:
$$
\|f(A)-f(B)\|\le\const\s\|f\|_{L^\be}\|A-B\|.
$$

To obtain those estimates techniques of double operator integrals and projective tensor products were used.

In this section we offer an elementary approach that shows that the above inequalities hold with constant equal to 1. Moreover, we obtain even sharper inequalities that can be considered as operator analogs of Bernstein's inequalities.

We also prove in this section that in general the operator moduli of continuity
$\O_f$ and $\O_f^\flat$ do not have to coincide.

\medskip

{\bf 5.1. Bernstein type inequalities for functions of self-adjoint operators.}
Let $\s>0$. Denote by $\mathscr E_\s$ the set of entire functions
of exponential type at most $\s$. The famous Bernstein theorem says that
$$
\sup_{x\in \R}|f^\prime(x)|\le\s\sup_{x\in \R}|f(x)|
$$
for every $f\in\mathscr E_\s$ (we refer the reader to \cite{L} for Bernstein's inequality, its generalizations
and related topics).
Bernstein's inequality implies that
\bay
\label{corber}
|f(x)-f(y)|\le\s\|f\|_{L^\infty}|x-y|
\ey
for every $f\in\mathscr E_\s$ and $x,y\in\R$, where $\|f\|_{L^\infty}\df\sup\limits_{x\in \R}|f(x)|$.

Bernstein \cite{B} also proved the following improvement of inequality \rf{corber}
\bay
\label{bern2}
\label{imprber}
|f(x)-f(y)|\le\b_\s(|x-y|)\|f\|_{L^\infty}
\ey
for every $f\in\mathscr E_\s$ and $x,y\in\R$,
where
$$
\b_\s(\d)=\left\{\begin{array}{ll}2\sin(\s \d/2),&\text {if}\,\,\,\,0\le \d\le\pi/\s,\\[.2cm]
2,&\text {if}\,\,\,\,\d>\pi/\s.
\end{array}\right.
$$
Iterating $m$ times inequality \rf{imprber}, we obtain
\bay
\label{bern3}
\label{iterber}
\|\D^m_hf\|_{L^\infty}\le\b_\s^m(|h|)\|f\|_{L^\infty}
\ey
for every $f\in\mathscr E_\s$ and $h\in\R$. This estimate is sharp,
because we have equality for $f(x)=e^{i\s x}$.

For the reader's convenience, we present the proof of Bernstein's inequality \rf{bern2}.

\medskip

{\bf Proof of \rf{bern2}.} It suffices to verify that
$|f(t)-f(-t)|\le2\sin t\|f\|_{L^\be}$ for all $f\in\mathscr E_1$ and
$t\in(0,\pi/2)$.  It is well known that the family $\left\{\dfrac{\cos z}{z-\frac\pi2-k\pi}\right\}_{k\in\Z}$
forms an orthogonal basis in the space $\mathscr E_1\cap L^2(\R)$, and
\bay
\label{Fz}
F(z)=\sum_{k\in\Z}(-1)^{k+1}F\Big(\frac\pi2+\pi k\Big)\frac{\cos z}{z-\frac\pi2-k\pi}
\ey
for all $F\in\mathscr E_1\cap L^2(\R)$, see, for example, \cite{L}, Lect. 20.2, Th. 1.
Applying \rf{Fz} to $F(z)=(f(z)-f(-z))z^{-1}$, we obtain
\bay
\label{vspo}
f(z)-f(-z)=\sum_{k\in\Z}\left(f\Big(\frac\pi2+\pi k\Big)-
f\Big(-\frac\pi2-\pi k\Big)\right)\frac{(-1)^{k}z\cos z}{(\frac\pi2+k\pi)(\frac\pi2+k\pi-z)},
\ey
whence
$$
|f(t)-f(-t)|\le2\|f\|_{L^\infty}\sum_{k\in\Z}\frac{t\cos
t}{(\frac\pi2+k\pi)(\frac\pi2+k\pi-t)}=2\|f\|_{L^\infty}\sin t.
$$
The last equality is an immediate consequence of \rf{vspo}. $\bl$

All inequalities stated at the beginning of this section are also true for entire functions with values
in a Banach space. To state the corresponding result, we need some notation.

Note that for every
$f\in\mathscr E_\s\cap L^\infty(\R)$,
\bay
\label{fEs}
|f(z)|\le e^{\s|\im z|}\|f\|_{L^\infty},\quad z\in\C,
\ey
see, for example, \cite{L}, page~97.
Let $\mu$ be a complex Borel measure on $\C$ such that $\int_\C e^{\s|\im z|}\,d|\mu|(z)<\infty$.
Then $\mu$ induces the following continuous linear functional:
$$
f\mapsto\int_\C f(z)\,d\mu(z)
$$
on $\mathscr E_\s\cap L^\infty(\R)$. Put
$$
\|\mu\|_{[\s]}\df\sup\left\{\left|\int_\C f\,d\mu\right|:f\in\mathscr E_\s, \|f\|_{L^\infty}\le1\right\}.
$$
Note that $\|\d_\l\|_{[\s]}=e^{\s|\im \l|}$, where $\d_\l$ denotes the $\d$-measure
at $\l$. The inequality $\|\d_\l\|_{[\s]}\le e^{\s|\im \l|}$ follows
from \rf{fEs}, while the opposite inequality is evident.
Inequality \rf{bern3} and the fact that it turns into equality for $f(x)=e^{i\s x}$ imply that
$$
\left\|\sum_{j=0}^{n}(-1)^{n-j}{n\choose j}\d_{t+jh}\right\|_{[\s]}=\b_\s^n(|h|)
$$
for all $t,h\in\R$.

Let $X$ be a complex Banach space. Denote by $\mathscr E_\s(X)$
the set of all $X$-valued
entire functions of exponential type at most $\s$. The definition of
entire functions of exponential type at most $\s$ with values in a Banach space is the same as in the case of scalar functions, see \cite{L}, Lect. 6.2.
Given $f\in \mathscr E_\s(X)$, we put
$$
\|f\|_{L^\infty}\df\sup\{\|f(t)\|_{X}:t\in\R\}.
$$

\begin{lem}
\label{0431}
Let $\s>0$ and let $\mu$ be a complex Borel measure on $\C$ such that
$$
\int_\C e^{\s|\im z|}\,d|\mu|(z)<\infty.
$$
Then for an arbitrary Banach space $X$ and for every
$f\in\mathscr E_\s(X)\cap L^\infty(\R,X)$,
\bay
\label{per}
\int_\C\|f(z)\|_{X}\,d|\mu|(z)<\be
\ey
and
\bay
\label{vto}
\left\|\int_\C f(z)\,d\mu(z)\right\|_{X}\le\|\mu\|_{[\s]}\|f\|_{L^\infty}.
\ey
\end{lem}

\Pf Let us first prove \rf{vto} under the assumption that \rf{per} holds. Indeed, it suffices to observe that
$$
\left\|\int_\C f(z)\,d\mu(z)\right\|_{X}
=\sup\left\{\left|\int_\C (f(z),u)\,d\mu(z)\right|:~u\in X^*,~\|u\|_{X^*}=1\right\}.
$$

Obviously, \rf{per} holds for measures $\mu$ with
compact support and, as we have just proved, \rf{vto} holds for such measures.

Applying \rf{vto} for $\mu=\d_\z$, where $\z\in\C$, we find that
\bay
\label{vne}
\|f(\z)\|_{X}\le e^{\s|\im \z|}\|f\|_{L^\infty}
\ey
for every $f\in\mathscr E_\s(X)\cap L^\infty(\R,X)$ and $\z\in\C$.
This immediately implies \rf{per}. $\bl$

\begin{lem}
\label{0434}
Let $A$ and $K$ be bounded self-adjoint operators. Suppose that $\s>0$ and
 $f\in\mathscr E_\s\cap L^\infty(\R)$. Then the operator-valued
function $\z\mapsto f(A+\z K)$ is an entire function of exponential type at most $\s\|K\|$.
\end{lem}

\Pf Put $\Phi(\z)=f(A+\z K)$. Let $\e>0$.
Applying von Neumann's inequality, we obtain
\bey
\|\Phi(\z)\|\le\max\{|f(z)|:|z|\le\|A\|+|\z|\|K\|\}
\le C_\e e^{(\s+\e)(\|A\|+|\z|\|K\|)}
\eey
for some constant $C_\e$.
Hence, $\Phi$ is an operator-valued entire function of exponential
type at most $\s\|K\|$. $\bl$

\begin{thm}
\label{phi}
Let $A$ and $K$ be self-adjoint operators with $\|K\|<\infty$. Suppose that
$\s>0$ and $f\in\mathscr E_\s\cap L^\infty(\R)$. Then there
exists an operator-valued entire function $\Phi$ of exponential type at most
$\s\|K\|$ such that $\Phi(t)=f(A+tK)$ for all $t\in \R$.
\end{thm}

\Pf Let $\{A_j\}_{j=1}^\infty$ be a sequence
of bounded self-adjoint operators such that
$
\lim\limits_{j\to\be}\|A_ju-Au\|=0
$
for every $u\in\cd_A$. Put $\Phi_j(\z)=f_j(A+\z K)$.
By Lemma \ref{0434}, $\Phi_j$ is an operator-valued function
of exponential type at most $\s\|K\|$. Moreover,
$\|\Phi_j\|_{L^\infty}\le\|f\|_{L^\infty}$ for every $j$.
By \rf{vne}, $\|\Phi_j(\z)\|\le e^{\s|\im\z|\cdot\|K\|}\|f\|_{L^\infty}$
for every $\z\in\C$ and $j\ge1$.
Note that $\lim\limits_{j\to\be}\Phi_j(t)=f(A+tK)$ in the strong operator topology
by Lemma 8.4 in \cite{AP2}.
Applying the vector version of the Vitali theorem (see, for e.g., \cite{HP}, Th. 3.14.1),
we find that
the sequence
$\{\Phi_j(\z)\}$ converges in the strong operator topology
for every $\z\in\C$ and the function $\Phi$ defined by
$\Phi(\z)\df\lim\limits_{j\to\be}\Phi_j(\z)$, $\z\in\C$, is an entire function.
It remains to observe that $\Phi(t)=f(A+tK)$ for all $t\in \R$
and $\|\Phi(\z)\|\le e^{\s|\im\z|\cdot\|K\|}\|f\|_{L^\infty}$
for every $\z\in\C$. $\bl$

\begin{thm}
\label{opbern}
Let $A$ and $B$ be self-adjoint operators such that $A-B$ is bounded. Then
$$
\|f(A)-f(B)\|\le\b_\s(\|A-B\|)\|f\|_{L^\infty}\le\s\|f\|_{L^\infty}\|A-B\|
$$
for every $f\in\mathscr E_\s$.
\end{thm}

\Pf By Theorem \ref{phi}, there exists  an operator-valued entire function $\Phi$ of exponential type at most
$\s\|B-A\|$
such that $\Phi(t)=f(A+t(B-A))$ for all $t\in \R$.
Inequality \rf{bern2} and Lemma \ref{0431} imply that
$$
\|f(A)-f(B)\|=\|\Phi(1)-\Phi(0)\|\le\b_{\s\|B-A\|}(1)\|\Phi\|_{L^\be}
\le\b_\s(\|A-B\|)\|f\|_{L^\be}. \quad\bl
$$

Recall that the inequality
$$
\|f(A)-f(B)\|\le\const \s\|f\|_{L^\infty}\|A-B\|
$$
for $f\in\mathscr E_\s$
is a special case of results of \cite{Pe3}.

\medskip

{\bf Remark.} It is natural to ask the question of whether the stronger inequality
\bay
\label{strongb}
\|f(A)-f(B)\|\le\const \|f^\prime\|_{L^\infty}\|A-B\|
\ey
holds for every $f\in\mathscr E_\s$. It turns out that the answer {\it is negative}.
Moreover, for every $\s>0$, there exists a function $f\in\mathscr E_\s$
such that $f^\prime$ is bounded on $\R$ and $\O_f(\d)=\be$ for all $\d\in(0,\be)$. 

\medskip

\Pf Let $\s>0$. Put
$$
f(x)\df x\int_0^{\s x}\frac{\sin t\,dt}t=x\int_0^\s\frac{\sin sx\,ds}s.
$$
It is easy to see that $f\in\mathscr E_\s$ and $f^{\prime}$ is bounded on $\R$.
It remains to note that $\O_f=\be$ by Theorem \ref{dif+-} or \ref{difbe}. $\bl$

Theorem \ref{usr} allows us to obtain the following consequence of
Theorem \ref{opbern}.

\begin{thm}
\label{opbern+}
Let $A$ and $B$ be self-adjoint operators and let $f\in\mathscr E_\s$. Then
for an arbitrary bounded operator $R$, the following inequality holds:
$$
\|f(A)R-Rf(B)\|\le\s\|f\|_{L^\infty}\|AR-RB\|.
$$
\end{thm}

\begin{thm}
\label{opbernsp}
Let $A$ and $B$ be self-adjoint operators and let $f\in\mathscr E_\s$. Suppose that $\fI$ is a quasinormed
ideal with majorization property. Then
for an arbitrary bounded operator $R$, the following inequality holds:
$$
\|f(A)R-Rf(B)\|_{\fI}\le\s\|f\|_{L^\infty}\|AR-RB\|_{\fI}.
$$
\end{thm}

\Pf The result follows from Theorem \ref{opbern+} and Theorem \ref{fI}. $\bl$

\begin{cor}
Let $A$ and $B$ be self-adjoint operators.  and let $f\in\mathscr E_\s$. Then
for an arbitrary bounded operator $R$, the following inequality holds:
$$
\|f(A)R-Rf(B)\|_{\bS_p}\le\s\|f\|_{L^\infty}\|AR-RB\|_{\bS_p},\quad 1\le p\le\be.
$$
\end{cor}

\begin{thm}
\label{dnest}
Let $A$ be a self-adjoint operator and $K$ a bounded self-adjoint operator  on Hilbert space.
Then for every $f\in\mathscr E_\s\cap L^\infty(\R)$, the following inequality holds:
$$
\|(\D^m_Kf)(A)\|\le\b_\s^m(\|K\|)\|f\|_{L^\infty}\le\s^m\|K\|^m\|f\|_{L^\infty}
$$
for every positive integer $m$.
\end{thm}

\Pf By Theorem \ref{phi}, there exists  an operator-valued entire function $\Phi$ of exponential type at most $\s\|K\|$
such that $\Phi(t)=f(A+tK)$ for every $t\in \R$.
Inequality \rf{iterber} and Lemma \ref{0431} imply that
$$
\|(\D^m_Kf)(A)\|=\|(\D^m_1\Phi)(0)\|\le\b_{\s\|B-A\|}^m(1)\|\Phi\|_{L^\be}
\le\b_\s^m(\|A-B\|)\|f\|_{L^\be}. \quad\bl
$$

\medskip

{\bf Remark.} Note that the inequality
$$
\|(\D^m_Kf)(A)\|\le\const\s^m\|K\|^m\|f\|_{L^\infty}
$$
is a consequence of results of \cite{Pe5} and \cite{AP2}.

\medskip

{\bf 5.2. Bernstein type inequalities and operator moduli of continuity.}
Let $f$ be a continuous function on $\R$. Denote by $\o_f$ the usual (scalar)
modulus of continuity of $f$,
$$
\o_f(\d)\df\sup\big\{|f(x)-f(y)|:~x,y\in\R,~\|x-y\|<\d\big\}.
$$
Clearly, $\o_f\le\O_f$. Note that $\O_f=\o_f$ for any affine function $f:\R\to\C$.
Theorem \ref{opbern} allows us to construct some more examples of such functions $f$.

\begin{thm}
\label{om0}
Let $f(t)=c_1 e^{{\rm i}\s t}+c_2 e^{-{\rm i}\s t}+c_3$, where $c_1,c_2,c_3\in\C$ and $\s>0$. Then
$\o_f=\O_f=(|c_1|+|c_2|)\b_\s$.
\end{thm}

\Pf We may assume that $\s=1$, $c_3=0$, and $|c_1|+|c_2|=1$. Note
that for every $h>0$, we have
\begin{align*}
\sup_t|f(t+h)-f(t)|&=\sup_t\big|c_1(e^{{\rm i} h}-1)e^{{\rm i} t}+c_2(e^{-{\rm i} h}-1)e^{-{\rm i} t}\big|\\[.2cm]
&=|c_1|\cdot|e^{{\rm i} h}-1|+|c_2|\cdot|e^{-{\rm i} h}-1|\\[.2cm]
&=2(|c_1|+|c_2|)|\sin(h/2)|=2|\sin(h/2)|.
\end{align*}
Hence, $\o_f(\d)=\sup\limits_{0<h\le \d}2|\sin(h/2)|=\b_1(\d)$.
Applying Theorem \ref{opbern}, we obtain
$$
\o_f(\d)\le\O_f(\d)\le(|c_1|+|c_2|)\,\b_1(\d)=\o_f(\d). \quad\bl
$$

\begin{thm}
\label{om3}
Let $f(t)=c_1 e^{{\rm i}\s t}+c_2 e^{-{\rm i}\s t}+c_3$, where $c_1,c_2,c_3\in\C$ and $\s>0$.
Then
$\O_f^{\flat}(\d)=(|c_1|+|c_2|)\min(2,\s \d)$.
\end{thm}

\Pf It suffices to consider the case where $\s=1$,
$c_3=0$,  and $|c_1|+|c_2|=1$. Moreover, taking into account that
$\O_f^{\flat}$ is invariant under translations of $f$, we may also
assume that $c_1,\,c_2\ge0$.
The inequality
$\O_f^{\flat}(\d)\le \d$ follows from Theorem \ref{opbern+},
and the inequality $\O_f^{\flat}\le 2$ is trivial.
It remains to prove that $\O_f^{\flat}(\d)\ge\min(2,\d)$.
It suffices to consider the case where $\d=2$, since the function $\O_f^{\flat}$ is nondecreasing and the function $t\mapsto t^{-1}\O_f^{\flat}(t)$ is nonincreasing.

Put $\f(t)\df 1-2\int_0^t\sgn\sin\pi s\,ds$.
We have $-1\le\f\le1$, $\f(t+1)=-\f(t)$ and $|\f^{\,\prime}|=2$ almost everywhere
on $\R$.
Let $M_\f$ be multiplication by $\f$ on $L^2(\R)$: $M_\f w\df\f w$.
Clearly, $\|M_\f\|=1$.
Let $A=-{\rm i}\dfrac{d}{dt}$ be defined on the set of functions
$h\in L^2(\R)$ such that $h^\prime\in L^2(\R)$.
Clearly, $AM_\f-M_\f A=-{\rm i}M_{\f^{\,\prime}}$. Hence, $\|AM_\f-M_\f A\|=2$,
whence $\O_f^{\flat}(2)\ge\|f(A)M_\f-M_\f f(A)\|$.

It is easy to see that $(e^{{\rm i}A}w)(t)=w(t+1)$ and $(e^{-{\rm i}A}w)(t)=w(t-1)$.
Hence, $f(A)w=c_1 w(t+1)+c_2 w(t-1)$. Thus
\begin{align*}
f(A)M_\f w&=c_1 \f(t+1)w(t+1)+c_2 \f(t-1)w(t-1)\\[.2cm]
&=-\f(t)\big(c_1 w(t+1)+c_2w(t-1)\big)=-M_\f f(A)w,
\end{align*}
whence $\|f(A)M_\f-M_\f f(A)\|=2\|M_\f f(A)\|$.

Let us show that $\|M_\f f(A)\|\ge1$.
Fix $\e\in(0,1/2)$. Denote by $P_\e$ the orthogonal projection onto
the space of functions $h\in L^2(\R)$ vanishing outside the $\e$-neighborhood
of $\Z$. Note that $P_\e(L^2(\R))$ is an invariant subspace of
the operators $M_\f$ and $f(A)$. We claim that $\|f(A)P_\e\|=1$
for every $\e>0$. Indeed, this follows from the following simple inequality:
$$
\lim_{n\to\be}\frac{\big\|f(A)P_\e\chi_{[0,n)}\big\|_{L^2}}{\big\|P_\e\chi_{[0,n)}\big\|_{L^2}}=1.
$$
It remains to observe that $\|M_\f w\|_{L^2}\ge(1-2\e)\|w\|_{L^2}$
for all $w\in P_\e(L^2(\R))$,
because $|\f|\ge1-2\e$ on the $\e$-neighborhood of $\Z$. $\bl$

\begin{cor}
\label{410}
Let $f(t)=e^{{\rm i}\s t}$, $\s>0$.
Then $\O_f(\d)<\O_f^{\flat}(\d)$
for all $\d\in\big(0,\pi\s^{-1}\big)$.
\end{cor}

\Pf It suffices to observe that $\O_f(\d)=2\sin(\s \d/2)$ for
$\d\in\big(0,\pi \s^{-1}\big)$
by Theorem \ref{om0} and $\O_f^{\flat}(\d)=\min(2,\s \d)$
for every $\d>0$ by Theorem \ref{om3}. $\bl$

\medskip

{\bf 5.3. Bernstein type inequalities in the case of unitary operators.}
We need the following elementary lemma.

\begin{lem}
\label{argexp}
Let $U$ and $V$ be unitary operators. Then there exists
a self-adjoint operator $A$ such that $V=e^{{\rm i}A}U$, $\|A\|\le\pi$ and $2\sin\big(\frac12\|A\|\big)=\|U-V\|$.
\end{lem}

\Pf Put $A=\arg (VU^{-1})$, where $\arg:\T\to\R$ is defined
by the formula $\arg(e^{is})=s$ for $s\in[-\pi,\pi)$.
Clearly, $2\sin\big(\frac12\|A\|\big)=\|I-VU^{-1}\|=\|U-V\|$. $\bl$

\begin{thm}
\label{uniber}
Let $f$ be a trigonometric polynomial of degree at most $d$.
Suppose that $|f(\z)|\le1$ for all $\z\in\T$.
Then
$$
\|f(U)-f(V)\|\le d\|U-V\|
$$
for every unitary operators $U$ and $V$.
\end{thm}

\Pf By Lemma \ref{argexp}, there exists a self-adjoint operator $A$
such that $V=e^{{\rm i}A}U$, $\|A\|\le\pi$ and
$2\sin\big(\frac12\|A\|\big)=\|U-V\|$.
Put $\Phi(\l)\df f\big(e^{{\rm i}\l A}U\big)$. Clearly, $\Phi$ is
an operator-valued entire function of exponential type at most $d\|A\|$
and $\|\Phi\|\le1$ on $\R$. Hence,
$$
\|f(U)-f(V)\|=\|\Phi(1)-\Phi(0)\|\le \b_{d\|A\|}(1)
$$
by inequality \rf{bern2} and Lemma \ref{0431}. It remains to observe that
$$
\b_{d\|A\|}(1)=\b_{\|A\|}(d)\le d \b_{\|A\|}(1)=2d\sin\left(\frac12\|A\|\right)=d\|U-V\|.\quad\bl
$$

Recall that the inequality
$$
\|f(U)-f(V)\|\le\const d\|U-V\|
$$
was proved in \cite{Pe1}.

\medskip

{\bf Remark.} Put
$$
\o_d(\d)\df\sup\big\{|z_1^d-z_2^d|:~z_1\in\T,~z_2\in\T,~
|z_1-z_2|<\d\big\}=\b_d\big(2\arcsin(\d/2)\big)
$$
for $\d\in(0,2]$. It is clear from the proof of Theorem \ref{uniber} that
$$
\|f(U)-f(V)\|\le\b_{d\|A\|}(1)=\b_d(\|A\|),
$$
where $A$ denote the same as in Theorem \ref{uniber}. Taking into account the fact that \lb$\|A\|=2\arcsin\big(\frac12\|U-V\|\big)$, we obtain the following sharp inequality
$$
\|f(U)-f(V)\|\le\b_d\big(2\arcsin(\|U-V\|/2)\big)=\o_d(\|U-V\|)
$$
under the hypotheses of Theorem \ref{uniber}.

\medskip

Theorem \ref{usru} allows us to obtain the following consequence of
Theorem \ref{uniber}.

\begin{thm}
\label{opbern+u}
Let $U$, $V$  and $f$ denote the same in Theorem {\em\ref{uniber}}. Then
$$
\|f(U)R-Rf(V)\|\le d\|UR-RV\|
$$
for every bounded operator $R$.
\end{thm}

\begin{thm}
\label{opbernspu}
Let $U$, $V$, $R$ and $f$ denote the same as in Theorem {\em\ref{opbern+u}}.
Then
$$
\|f(U)R-Rf(V)\|_{\fI}\le d\|UR-RV\|_{\fI}
$$
for every quasinormed ideal $\fI$ with majorization property.
\end{thm}
\Pf The result follows from Theorem \ref{opbern+u} and Theorem \ref{fIu}. $\bl$

\begin{cor}
Let $U$, $V$, $R$ and $f$ denote the same as in Theorem {\em\ref{opbern+u}}.
Then
$$
\|f(U)R-Rf(V)\|_{\bS_p}\le d\|UR-RV\|_{\bS_p}
$$
for every $p\in[1,\be]$.
\end{cor}



\begin{thm}
Let $A$ be a self-adjoint operator and let $U$ be a unitary operator on Hilbert space. Then
$$
\left\|\sum_{k=0}^m(-1)^k\left(\begin{matrix}m\\k\end{matrix}\right)f\big(e^{{\rm i}kA}U\big)
\right\|\le d^m\|e^{{\rm i}A}-I\|^m\le d^m\|A\|^m.
$$
for every positive integer $m$ and every trigonometric polynomial $f$ of degree at most $d$ with $|f|\le1$
on $\T$.
\end{thm}

\Pf By Lemma \ref{argexp}, we can assume that
$\|e^{{\rm i}A}-I\|=2\sin(\frac12\|A\|)$.
Put $\Phi(\l)\df f(e^{{\rm i}\l A}U)$. Clearly, $\Phi$ is
an operator-valued entire function of exponential type at most $d\|A\|$
and $\|\Phi\|\le1$ on $\R$. Hence,
\begin{align*}
\left\|\sum_{k=0}^m(-1)^k\left(\begin{matrix}m\\k\end{matrix}\right)f\big(e^{{\rm i}kA}U\big)
\right\|&=\|\D_1^m\Phi(0)\|\le\b_{d\|A\|}^m(1)\\[.2cm]
&=\b_{\|A\|}^m(d)\le d^m\b_{\|A\|}^m(1)
=d^m\|e^{{\rm i}A}-I\|^m
\end{align*}
by inequality \rf{iterber} and Lemma \ref{0431}. $\bl$

\medskip

{\bf Remark.} The inequalities
$$
\left\|\sum_{k=0}^m(-1)^k\left(\begin{matrix}m\\k\end{matrix}\right)f\big(e^{{\rm i}kA}U\big)
\right\|\le\const d^m\|e^{{\rm i}A}-I\|^m\le\const d^m\|A\|^m
$$
are consequences of results of \cite{Pe5} and \cite{AP2}.

\

\section{\bf H\"older--Zygmund estimates for self-adjoint operators}
\setcounter{equation}{0}
\label{s5}

\

In this section we show that the estimates of operator finite differences that were obtain in \cite{AP2} for bounded self-adjoint operators also hold in the case of arbitrary (not necessarily bounded) self-adjoint operators.

The following theorem was proved in \cite{AP2}.

\begin{thm}
\label{saH}
Let $0<\a<1$. Then there is a constant $c>0$ such that for
every $f\in\L_\a(\R)$ and for arbitrary self-adjoint operators $A$ and $B$ on Hilbert space
the following inequality holds:
$$
\|f(A)-f(B)\|\le c\,\|f\|_{\L_\a(\R)}\|A-B\|^\a.
$$
\end{thm}

To state the result for arbitrary H\"older--Zygmund classes $\L_\a(\R)$ with $\a>0$, we introduce the following notation
\bay
\label{defD}
\big(\D_K^mf\big)(A)\df\sum_{j=0}^m(-1)^{m-j}\left(\begin{matrix}m\\j\end{matrix}\right)f\big(A+jK\big)
\ey
for functions $f\in C(\R)$ at least in the case when $A$ and $K$ are bounded
self-adjoint operators. It is also clear that formula \rf{defD} can be used for bounded functions $f$ in $\L_\a(\R)$ even if $A$ is unbounded.

However, for $f\in\L_\a(\R)$ with $\a\ge1$, in the case of unbounded self-adjoint operators $A$ it may happen that
$\cd_{f(A)}\cap\cd_{f(A+K)}=\{0\}$, see Corollary \ref{emptya} below. Nevertheless, by Theorem \ref{saH}, for $f\in\L_\a(\R)$ with $\a<1$,
we have $\cd_{f(A)}=\cd_{f(A+K)}$ . Thus if we want to consider
finite differences $\big(\D_K^mf\big)(A)$ also for unbounded self-adjoint operators $A$
in the case $f\in\L_\a(\R)$ with $\a\ge1$, we have to
define $\big(\D_K^mf\big)(A)$ more accurately.

Let us first show for every $\a\ge1$, there exists $f\in\L_\a(\R)$ such that the intersection $\cd_{f(A)}\cap\cd_{f(A+K)}$ is trivial.

\begin{lem}
\label{domA}
Let $A$ be an unbounded self-adjoint operator. Then there
exists an orthogonal projection $P$ such that
$\cd_A\cap P^{-1}(\cd_A)=\{0\}$.
\end{lem}

\Pf First we consider the following special case where
$A_0\f\df-{\rm i}\f'$ in the space $L^2[0,2\pi]$ with domain
$$
\cd_{A_0}=\big\{\f\in L^2[0,2\pi]:~\f'\in L^2[0,2\pi],~\f(0)=\f(2\pi)\big\}.
$$
Let $E$ be a Lebesgue measurable subset of $[0,2\pi]$
such that each of the sets $\D\cap E$ and $\D\setminus E$ has positive Lebesgue measure for every nondegenerate interval $\D$ in $[0,2\pi]$.
Let $P$ denote multiplication by the characteristic function $\chi_E$ of $E$.
Clearly, $P$ is an orthogonal projection and $\cd_{A_0}\cap P^{-1}(\cd_{A_0})=\{0\}$, because $\cd_{A_0}\subset C([0,2\pi])$.

Let now $A=f(A_0)$ where $f$ is a real continuous on $\R$ function such that
\lb $\lim\limits_{|t|\to\be}t^{-1}|f(t)|=\be$. Then $\cd_A\cap P^{-1}(\cd_A)=\{0\}$
because $\cd_A\subset\cd_{A_0}$.

If $A$ is an arbitrary unbounded self-adjoint operator with pure point spectrum, then it is unitarily equivalent to an orthogonal sum
$A_\be=\bigoplus\limits_{j=1}^\be f_j(A_0)$ of operators considered above.
Clearly, $\cd_{A_\be}\cap P_\be^{-1}(\cd_{A_\be})=\{0\}$ for $P_\be=\bigoplus\limits_{j=1}^\be P$.

Finally, it remains to observe that for every unbounded self-adjoint operator $A$
there exists a self-adjoint operator $A_\flat$ with pure point spectrum
such that $A-A_\flat$ is bounded, and so $\cd_A=\cd_{A_\flat}$. $\bl$

It is easy to see that if $f\in C(\R)$ and $|f(t)|\le\const(1+|t|)$, then $\cd_{f(A)}\supset\cd_{A}$.
Hence, $\cd_{f(A)}\cap\cd_{f(A+K)}$ is a dense subset for every bounded
self-adjoint operator $K$.

\begin{thm}
\label{empty}
Let $f$ be a real function continuous on $\R$. Suppose that
\lb $\lim\limits_{t\to+\be}t^{-1}|f(t)|=\be$. Then there exists a self-adjoint
operator $A$ and an orthogonal projection $P$ such that
$\cd_{f(A)}\cap\cd_{f(A+P)}=\{0\}$.
\end{thm}

\Pf Note that $\cd_{f(A)}=\cd_{|f|(A)}$. Hence, we may assume that $f\ge0$.
Let us first consider the special case when there exists a function
$g\in\L_{1/2}(\R)$ such that $f(t)=tg(t)$ for $t\ge1$.
Let $A$ be a self-adjoint operator with $\s(A)=[1,\be)$.
By Lemma \ref{domA},
there exists an orthogonal projection $P$ such that $\cd_{g(A)}\cap P^{-1}(\cd_{g(A)})=\{0\}$.
Clearly, $\s(A+P)\subset[1,\be)$.
Let us prove that $\cd_{f(A)}\cap\cd_{f(A+P)}=\{0\}$.
Suppose that $u\in\cd_{f(A)}\cap\cd_{f(A+P)}$. Then $Au\in\cd_{g(A)}$ and
$Au+Pu\in\cd_{g(A+P)}$. Note that $\cd_{g(A)}=\cd_{g(A+P)}$
by Theorem \ref{saH}. Hence, $Pu\in\cd_{g(A)}$, while
$u\in\cd_{f(A)}\subset\cd_{g(A)}$, and so $u=0$ because of the equality
$\cd_{g(A)}\cap P^{-1}(\cd_{g(A)})=\{0\}$.

To complete the proof, it suffices verify that there exists a function
$g\in\L_{1/2}(\R)$ such that $\lim\limits_{t\to\be}g(t)=\be$
and $f(t)\ge tg(t)$. We can assume that the function $t\mapsto t^{-1}f(t)$
is nondecreasing on $[1,\be)$. Let $\f$ be a nondecreasing function in $\L_{1/2}(\R)$ such that
$\f(t)=0$ for $t\le0$ and $\f(t)=1$ for $t\ge1$. For $a,b\in\R$ with $a<b$, we put
$\f_{a,b}(t)\df\f\big(\frac{t-a}{b-a}\big)$.
Clearly, $\|\f_{a,b}\|_{\L_{1/2}(\R)}=(b-a)^{-1/2}\|\f\|_{\L_{1/2}(\R)}$.
We can construct by induction a nondecreasing sequence $\{a_k\}_{k\ge0}$ of numbers
such that $a_0=1$ and
$$
\sum_{k=1}^\be\frac{a_k^{-1}f(a_k)-a_{k-1}^{-1}f(a_{k-1})}{(a_{k+1}-a_k)^{1/2}}<\be.
$$
Put
$$
g(t)=a_0^{-1}f(a_0)\f_{a_0,a_1}+\sum_{k=1}^\be\Big(a_k^{-1}f(a_k)-a_{k-1}^{-1}f(a_{k-1})\Big)\f_{a_k,a_{k+1}}.
$$
Clearly, $g\in\L_{1/2}(\R)$, $\lim\limits_{t\to\be}g(t)=\be$, and $tg(t)\le f(t)$ for $t\in[1,\be)$. $\bl$

\begin{cor}
\label{emptya}
Let $\a\ge1$. Then there exists a function $f\in\L_\a(\R)$, a self-adjoint
operator $A$, and an orthogonal projection $P$ such that $\cd_{f(A)}\cap\cd_{f(A+P)}=\{0\}$.
\end{cor}
\Pf If $\a>1$, we can apply Theorem \ref{empty} to $f(t)=|t|^\a$.
If $\a=1$, we can apply Theorem \ref{empty} to $f(t)=t\log|t|$. $\bl$

Now we return to the definition of $\big(\D_K^mf\big)(A)$.
We have already mentioned that formula \rf{defD} can be used in the case when $A$ is a bounded operator as well as in the case when $f$ is a bounded function.
Moreover, as we have observed above, if $f$ satisfies the inequality $|f(t)|\le\const(1+|t|)$, then
$\cd_{f(A+jK)}\supset\cd_A$ for all $j$ which allows us to define
the operator $\big(\D_K^mf\big)(A)$ by formula \rf{defD} on a dense subset.

For $f\in C(\R)\cap\mathscr S^\prime(\R)$, we will also use the following formula to define the finite difference $\big(\D_K^mf\big)(A)$:
\bay
\label{defD2}
\big(\D_K^mf\big)(A)\df\sum_{n\in\Z}\left(\sum_{j=0}^m(-1)^{m-j}
\left(\begin{matrix}m\\j\end{matrix}\right)f_n\big(A+jK\big)\right),
\quad f_n\df f*W_n+f*W_n^\sharp,
\ey
under the assumption that
$$
\sum_{n\in\Z}\left\|\sum_{j=0}^m(-1)^{m-j}
\left(\begin{matrix}m\\j\end{matrix}\right)f_n\big(A+jK\big)\right\|<\be.
$$

It will be clear from the Remark following Theorem \ref{fs} and from Theorem \ref{As} that this definition does not depend on the choice of the functions $W_n$.

\begin{thm}
\label{sam}
Let $0<\a<m$. Then there exists a constant $c>0$ such that for every self-adjoint operators  $A$ and $K$
with $\|K\|<\be$, the following inequality holds:
$$
\sum_{n\in\Z}\left\|\sum_{j=0}^m(-1)^{m-j}
\left(\begin{matrix}m\\j\end{matrix}\right)f_n\big(A+jK\big)\right\|\le c\,\|f\|_{\L_\a(\R)}\|K\|^\a,
$$
and so $\big(\D_K^mf\big)(A)$ is well defined by  {\rm\rf{defD2}} and
$$
\left\|\big(\D_K^mf\big)(A)\right\|
\le c\,\|f\|_{\L_\a(\R)}\cdot\|K\|^\a.
$$
\end{thm}

\Pf Theorem \ref{dnest} implies the following inequality
\begin{align*}
\left\|\sum_{j=0}^m(-1)^{m-j}\left(\begin{matrix}m\\j\end{matrix}\right)f_n\big(A+jK\big)
\right\|&\le2^{mn+m}\|K\|^m\|f_n\|_{L^\infty}\\[.2cm]
&\le\const2^{(m-\a)n}\|K\|^m\|f\|_{\L_\a(\R)}.
\end{align*}
Moreover, it is clear that
\begin{align*}
\left\|\sum_{j=0}^m(-1)^{m-j}\left(\begin{matrix}m\\j\end{matrix}\right)f_n\big(A+jK\big)
\right\|&\le\sum_{j=0}^m\left(\begin{matrix}m\\j\end{matrix}\right)\|f_n\big(A+jK\big)\|\\[.2cm]
&\le2^m\|f_n\|_{L^\infty}\le\const2^{-\a n}\|f\|_{\L_\a(\R)}.
\end{align*}
It remains to observe that
$$
\sum_{n\in\Z}2^{-\a n}\min(2^{mn}\|K\|^m,1)=\sum_{2^n\|K\|\le1}2^{(m-\a) n}\|K\|^m+
\sum_{2^n\|K\|>1}2^{-\a n}\le\const\|K\|^\a.\quad\bl
$$

\begin{thm}
\label{fs}
Let $0<\a<m$. Suppose that $\{f_j\}_{j=1}^\infty$ is a bounded sequence
of functions in $\L_\a(\R)$ that converges pointwise to a function $f$.
Then $f\in\L_\a(\R)$ and for every self-adjoint operators  $A$ and $K$ with
$\|K\|<\be$,
$$
\big(\D_K^mf\big)(A)=\lim\limits_{j\to\be}\big(\D_K^mf_j\big)(A)
$$
in the strong operator topology.
\end{thm}

\Pf By Theorem \ref{sam}, it suffices to verify that
$$
\lim_{s\to\be}(f_s*W_n+f_s*W_n^\sharp)(A+jK)=(f*W_n+f*W_n^\sharp)(A+jK)
$$
for all $n,j\in\Z$ in the strong operator topology. This follows
from the fact that 
$$
\sup_s\|f_s*W_n+f_s*W_n^\sharp\|_{L^\infty}<+\infty
$$
and from the fact that
$$
\lim\limits_{s\to\be}(f_s*W_n+f_s*W_n^\sharp)(t)=(f*W_n+f*W_n^\sharp)(t)
$$
for every $t\in\R$. $\bl$

\medskip

{\bf Remark.} This theorem allows us to give another (equivalent) definition
of the finite difference $\big(\D_K^mf\big)(A)$ for $f\in\L_\a(\R)$.
Let $f\in\L_a(\R)$. Then there exists a sequence $\{f_s\}_{s=1}^\infty$ in $\L_\a(\R)\cap L^\be(\R)$
such that $\sup\limits_s\|f_s\|_{\L_\a(\R)}<\be$ and $\lim\limits_{s\to\be}f_s(t)=f(t)$
for all $t\in\R$. Put
$$
\big(\D_K^mf\big)(A)\df\lim_{s\to\be}\left(\sum_{j=0}^m(-1)^{m-j}\left(\begin{matrix}
m\\j\end{matrix}\right)f_s\big(A+jK\big)\right),
$$
where the limit is taken in the strong operator topology.

\medskip

The following theorem yields one more (equivalent) definition
of $\big(\D_K^mf\big)(A)$ for \lb$f\in\L_\a(\R)$.

\begin{thm}
\label{As}
Let $0<\a<m$. Suppose that $A$ and $K$ are self-adjoint operators
such that $\|K\|<\be$. Let $\{A_s\}_{s=1}^\infty$ be a sequence of bounded
self-adjoint operators such that $\lim\limits_{s\to\be}\|A_su-Au\|=0$
for all $u\in\cd_A$. Then for $f\in\L_\a(\R)$,
$$
\big(\D_K^mf\big)(A)=\lim_{s\to\be}\left(\sum_{j=0}^m(-1)^{m-j}\left(\begin{matrix}
m\\j\end{matrix}\right)f\big(A_s+jK\big)\right),
$$
where the limit is taken in the strong operator topology.
\end{thm}

\Pf By Theorem \ref{sam}, it suffices to verify that
$$
\lim_{s\to\be}f_n(A_s+jK)=f_n(A+jK).
$$
in the strong operator topology,
where $f_n=f*W_n+f*W_n^\sharp$, $n,\,j\in\Z$.
This follows from Lemma 8.4 in \cite{AP2}.
$\bl$

\begin{thm}
\label{1355}
Let $0<\a<m$ and let $f$ be a functions in $\L_\a(\R)$
satisfying $|f(t)|\le\const(1+|t|)$.
Suppose that $A$ and $K$ are self-adjoint operators such that $\|K\|<\be$.
Then $\cd_{f(A+jK)}\supset\cd_A$ for every $j$. Moreover, if
$\big(\D_K^mf\big)(A)$ is defined by {\rm\rf{defD2}}, then
$$
\big(\D_K^mf\big)(A)u=\sum_{j=0}^m(-1)^{m-j}\left(\begin{matrix}m\\j\end{matrix}\right)f\big(A+jK\big)u
$$
for every $u\in\cd_A$.
\end{thm}

\Pf We can take a sequence $\{A_s\}_{s=1}^\infty$ of bounded
self-adjoint operators such that $\lim\limits_{s\to\be}\|A_su-Au\|=0$
for all $u\in\cd_A$. It remains to apply Theorem \ref{As} and Lemma 8.5 in \cite{AP2}. $\bl$


\medskip

{\bf Remark.} We can interpret Theorem \ref{sam} in the following way. Consider the measure $\nu$ on $\R$ defined by
$$
\nu\df\D_1^m\d_0=\sum_{j=0}^m(-1)^{m-j}\left(\begin{matrix}m\\j\end{matrix}\right)\d_{-j},
$$
where for $a\in\R$, $\d_a$ is the unit point mass at $a$. Then
$$
\sum_{j=0}^m(-1)^{m-j}\left(\begin{matrix}m\\j\end{matrix}\right)f\big(A+jK\big)=
\int_\R f(A-tK)\,d\nu(t)
$$
at least formally.
Clearly, $\nu$ determines a continuous linear functional on $\l_\a(\R)$ by the formula
$$
f\mapsto \int_\R f(t)\,d\nu(t).
$$
In other words, $\nu\in B_1^{-\a}(\R)$ (see \S\,\ref{s2}).
In \cite{AP2} we generalized Theorem \ref{sam} to the case of an arbitrary distribution $\nu$ in $B_1^{-\a}(\R)$ in the case of bounded $A$.

\medskip

Here we consider the case of an arbitrary self-adjoint operator $A$.
Let $a\in\R$ and $h>0$. Put $\nu_{a,h}\df\D^m_h\d_a$. Then we have formally
$$
\int_\R f(A-tK)\,d\nu_{a,h}(t)=\sum_{j=0}^m(-1)^{m-j}
\left(\begin{matrix}m\\j\end{matrix}\right)f\big(A-aK+jhK\big).
$$
Put
\bay
\label{defQ}
\Q_{A,K}^{\nu_{a,h}}f\df\sum_{j=0}^m(-1)^{m-j}
\left(\begin{matrix}m\\j\end{matrix}\right)f\big(A-aK+jhK\big).
\ey
Clearly,
$$
\|\Q_{A,K}^{\nu_{a,h}}f\|
\le c\,\|f\|_{\L_\a(\R)}h^\a\|K\|^\a=\const\|f\|_{\L_\a(\R)}\|\nu_{a,h}\|_{B_1^{-\a}(\R)}\|K\|^\a.
$$

Denote by $\frak L_m$ the linear span of the family $\{\D_h^m\d_a\}_{a,h}$.
Clearly, we can extend the definition of the operator $\Q_{A,K}^gf$
to the case when $g\in\frak L_m$ so that the operator
\lb$g\mapsto\Q_{A,K}^gf$ is linear.

\begin{thm}
\label{1366}
Let $\a>0$ and $m-1\le\a<m$ for $m\in\Z$. Suppose that $A$ and $K$ are self-adjoint operators such that $\|K\|<\be$.
Then the operator $g\mapsto\Q_{A,K}^gf$ initially defined on $\frak L_m$
by {\em\rf{defQ}}
extends to a continuous linear operator from $B_1^{-\a}(\R)$ to ${\mathscr B}(\h)$.
Moreover,
$$
\big\|\Q_{A,K}^gf\|\le\const\|f\|_{\L_\a(\R)}\|g\|_{B_1^{-\a}(\R)}\|K\|^\a ,
$$
where the constant depends only on $\a$.
\end{thm}

\Pf
We use the following result on atomic decomposition of $B_1^{-\a}(\R)$
(see \cite{A}, Ch. 3, Th. 3.1): if $g\in B_1^{-\a}(\R)$, then $g$ expands in a norm convergent series
$$
g=\sum_{j=1}^\be\l_j\D^m_{h_j}\d_{a_j},\quad \,a_j\in\R, ~h_j>0,
$$
such that
$$
\sum_{j=1}^\be|\l_j|h_j^\a\le\const\|g\|_{B_1^{-\a}(\R)}.
$$
Put
\bay
\label{korr}
\Q_{A,K}^gf\df\sum_{j=1}^\infty \Q_{A,K}^{g_j}f,
\ey
where $g_j\df\l_j\d_{a_j}$. Then
$$
\|\Q_{A,K}^gf\|\le\sum_{j=1}^\infty \|\Q_{A,K}^{g_j}f\|
\le c\,\|f\|_{\L_\a(\R)}\|K\|^\a\sum_{j=1}^\be |\l_j|h_j^\a\le C\,\|f\|_{\L_\a(\R)}\|g\|_{B_1^{-\a}(\R)}\|K\|^\a.
$$
Let us verify that $\Q_{A,K}^gf$ is well defined by \rf{korr}.
Suppose that
$$
\sum_{j=1}^\be|\l_j|h_j^\a<\be\quad \text {and}\quad\sum_{j=1}^\be\l_j\D^m_{h_j}\d_{a_j}=\0
\quad \text {in}\quad B_1^{-\a}(\R).
$$
Let us prove that $\sum\limits_{j=1}^\be \Q_{A,K}^{g_j}f=0$ in the space ${\mathscr B}(\h)$,
where $g_j=\l_j\d_{a_j}$.
By Theorem \ref{sam}, we have $\Q_{A,K}^{g_j}f=\sum\limits_{n\in\Z}\Q_{A,K}^{g_j}f_n$
with $\sum\limits_{n\in\Z}\|\Q_{A,K}^{g_j}f_n\|\le C\,|\l_j|h_j^\a\,\|f\|_{\L_\a(\R)}\|K\|^\a$.
In particular,
$$
\sum\limits_{j=1}^\be\sum_{n\in\Z}\|\Q_{A,K}^{g_j}f_n\|
\le C\,\|f\|_{\L_\a(\R)}\|K\|^\a\sum\limits_{j=1}^\be|\l_j|h_j^\a<\be.
$$
Hence,
\bey
\sum\limits_{j=1}^\be \Q_{A,K}^{g_j}f=\sum\limits_{j=1}^\be\sum_{n\in\Z} \Q_{A,K}^{g_j}f_n
=\sum_{n\in\Z}\sum\limits_{j=1}^\be \Q_{A,K}^{g_j}f_n.
\eey
It remains to verify that $\sum\limits_{j=1}^\be \Q_{A,K}^{g_j}f_n=0$ for all $n\in\Z$.
Let $u,v\in\h$.
Put $\f_n(t)\df (f_n(A-tK)u,v)$.
Clearly,
$$
\sum_{j=1}^\be (\Q_{A,K}^{g_j}f_nu,v)=\sum_{j=1}^\be\int_\R\f_n(t)dg_j(t)=0,
$$
because $f_n\in\L_\a(\R)$ by Theorem \ref{sam} and $\sum\limits_{j=1}^\be g_j=0$
in the space $B_1^{-\a}(\R)$. $\bl$

The following theorem is a generalization of Theorem \ref{fs}.
\begin{thm}
\label{fs+}
Let $g\in B_1^{-\a}(\R)$ with $\a>0$. Let $\{f_s\}_{s=1}^\infty$ be a bounded sequence
of function in $\L_\a(\R)$ which converges pointwise to a function $f$.
Then $f\in\L_\a(\R)$ and $\Q_{A,K}^gf=\lim\limits_{s\to\be}\Q_{A,K}^gf_s$
in the strong operator topology
for every self-adjoint operators  $A$ and $K$ with $\|K\|<\be$.
\end{thm}

\Pf Let $m$ be the smallest integer greater than $\a$. The case when $g=\D_h^m\d_a$
follows from Theorem \ref{fs}. In the general case we can use the atomic
decomposition of $B_1^{-\a}(\R)$, see the proof of Theorem \ref{1366}. $\bl$

In the same way we can obtain the following generalization of Theorem \ref{As}.

\begin{thm}
Let $A$ and $K$ be self-adjoint operators
such that $\|K\|<\be$. Suppose that $\{A_s\}_{s=1}^\infty$ is a sequence of bounded
self-adjoint operators such that $\lim\limits_{s\to\be}\|A_su-Au\|=0$
for $u\in\cd_A$. Then for every $\a>0$ and $f\in\L_\a(\R)$,
$$
\Q_{A,K}^gf=\lim\limits_{s\to\be}\Q_{A_s,K}^gf,
$$
where the limit is taken in the strong operator topology.
\end{thm}

\medskip

{\bf Remark 1.} All results of \S\,5 in \cite{AP3} are true also for not necessarily bounded
self-adjoint operators $A$ (and $B$). The corresponding finite differences
can be defined by \rf{defD2}, $K\df B-A$ in the case $m=1$.
Finite differences $(\D^m_Kf)(A)$
can also be defined with the help of approximation $f_s\to f$
or $A_s\to A$ as in the remark following Theorem \ref{fs} and in Theorem \ref{As}.

\medskip

{\bf Remark 2.} To obtain the results of \S\,11 in \cite{AP3} for not necessarily bounded
self-adjoint operators $A$ and $B$, we can apply
Lemma \ref{anbnrn}.


\

\section{\bf Higher order moduli of continuity}
\setcounter{equation}{0}
\label{s6}

\

In \cite{AP2} we obtained estimates of operator finite differences for functions
of class $\L_{\o,m}$ (see \S\,2) in the case of bounded self-adjoint operators. The purpose of this section is to extend the results of \cite{AP2} to the case of arbitrary self-adjoint operators.

Let
$\o$ be a nondecreasing function on $(0,\be)$ such that
\bay
\label{om}
\lim_{x\to0}\o(x)=0\quad\mbox{and}\quad
\o(2x)\le2^m\o(x)\quad\mbox{for}\quad x>0.
\ey

Denote by $\L_{\o,m}(\R)$ is the space of continuous functions $f$ on $\R$
satisfying
$$
\|f\|_{\L_{\o,m}(\R)}\df\sup\limits_{t>0}\frac{\|\D^m_tf\|_{L^\infty}}{\o(t)}<+\infty.
$$

In this section we obtain norm estimates for finite differences
$$
\big(\D_K^mf\big)(A)\df\sum_{j=0}^m(-1)^{m-j}\left(\begin{matrix}m\\j\end{matrix}\right)f\big(A+jK\big)
$$
for functions $f\in\L_{\o,m}(\R)$ and self-adjoint operators $A$ and $K$ with $\|K\|<\be$.
As mentioned in \S\,\ref{s5}, this definition of $\big(\D_K^mf\big)(A)$ does not work
in general
for unbounded $A$ and $m\ge2$. We will apply the following formula
\bay
\label{defD4}
\big(\D_K^mf\big)(A)\df\sum_{n=-\be}^N\big(\D_K^mf_n\big)(A)+\big(\D_K^m(f-f*V_N)\big)(A),
\ey
where $f_n\df f*W_n+f*W_n^\sharp$ and $V_N$ is defined in \S\,\ref{s2}. We say that \rf{defD4} well defines
the finite difference $\big(\D_K^mf\big)(A)$ if
$$
\sum_{n=-\be}^N\|\big(\D_K^mf_n\big)(A)\|<\be.
$$
Note that $f-f*V_N\in L^\infty(\R)$ and this definition of $\big(\D_K^mf\big)(A)$
does not depend on the choice of $N\in\Z$.

Given a nondecreasing function $\o$ satisfying \rf{om}, we define the function $\o_{*,m}$ by
$$
\o_{*,m}(x)=x^m\int_x^\be\frac{\o(t)}{t^{m+1}}\,dt=\int_1^\be\frac{\o(sx)}{s^{m+1}}\,dx.
$$

The following theorem in the case when $A$ is bounded was proved in \cite{AP2}.

\begin{thm}
\label{oon}
Let $m$ be a positive integer.  Then there exists a positive number $c$ such that for an arbitrary
nondecreasing function $\o$ on $(0,\be)$ satisfying {\em\rf{om}} and
such that $\o_{*,m}(x)<\be$ for $x>0$, an arbitrary
bounded function $f$ in $\L_{\o,m}(\R)$, and arbitrary self-adjoint operators $A$ and
$K$ with $\|K\|<\be$ the following inequality holds:
$$
\left\|\big(\D_K^mf\big)(A)\right\|\le c\,\|f\|_{\L_{\o,m}(\R)}\,\o_{*,m}\big(\|K\|\big),
$$
where the finite difference $\big(\D_K^mf\big)(A)$ is well defined by  {\rm\rf{defD4}}.
Moreover,
$$
\sum_{n=-\infty}^N\big\|\big(\D_K^mf_n\big)(A)\big\|+\big\|\big(\D_K^m(f-f*V_N)\big)(A)\big\|
\le C\,\|f\|_{\L_{\o,m}(\R)}\,\o_{*,m}\big(\|K\|\big)
$$
provided $2^{-N}<\|K\|\le2^{-N+1}$.
\end{thm}

\Pf
%
By Theorem \ref{mnn},
\begin{align*}
\big\|\big(\D_K^m(f-f*V_N)\big)(A)\big\|&\le\const\|f-f*V_N\|_{L^\be}\\[.2cm]
&\le\const\|f\|_{\L_{\o,m}(\R)}\o\big(2^{-N}\big)
\le\const\|f\|_{\L_{\o,m}(\R)}\o_{*,m}\big(||K\|\big).
\end{align*}
On the other hand, it follows from Theorem \ref{dnest}, and Corollary \ref{Wnm} that
$$
\big\|\big(\D_K^mf_n\big)(A)\big\|\le\const2^{mn}\|f_n\|_{L^\be}\|K\|^m
\le\const\|f\|_{\L_{\o,m}(\R)}2^{mn}\o\big(2^{-n}\big)\|K\|^m.
$$
Thus
\begin{align*}
\sum_{n=-\be}^N\big\|\big(\D_K^mf_n\big)(A)\big\|&\le\const\sum_{n=-\be}^N\|f\|_{\L_{\o,m}
(\R)}2^{mn}\o\big(2^{-n}\big)\|K\|^m\\[.2cm]
&=\sum_{k\ge0}2^{(N-k)m}\o\big(2^{N-k}\big)\|f\|_{\L_{\o,m}(\R)}\|K\|^m\\[.2cm]
&\le\const\left(\int_{2^{-N}}^\be\frac{\o(t)}{t^{m+1}}\,dt\right)\|f\|_{\L_{\o,m}(\R)}\|K\|^m\\[.2cm]
&=\const2^{-Nm}\o_{*,m}\big(2^{-N}\big)\|f\|_{\L_{\o,m}(\R)}\|K\|^m\\[.2cm]
&\le\const\|f\|_{\L_{\o,m}(\R)}\o_{*,m}\big(\|K\|\big).
\end{align*}
This completes the proof. $\bl$

\begin{thm}
\label{fs2}
Let $m$ and $\o$ satisfy the hypotheses of Theorem {\rm\ref{oon}}. Suppose that $\{f_s\}_{s=1}^\infty$ is a bounded sequence
of functions in $\L_{\o,m}(\R)$ that converges pointwise to a function $f$.
Then $f\in\L_{\o,m}(\R)$ and
for every self-adjoint operators  $A$ and $K$ with $\|K\|<\be$,
$$
\big(\D_K^mf\big)(A)=\lim\limits_{s\to\be}\big(\D_K^mf_s\big)(A)
$$
in the strong operator topology.
\end{thm}

The proof is similar to the proof of Theorem \ref{fs}.

%
%

\begin{thm}
\label{As2}
Let $m$ and $\o$ satisfy the hypotheses of Theorem {\rm\ref{oon}}. Let $A$ and $K$ be self-adjoint operators
with $\|K\|<\be$. Let $\{A_s\}_{s=1}^\infty$ be a sequence of bounded
self-adjoint operators such that $\lim\limits_{s\to\be}\|A_su-Au\|=0$
for all $u\in\cd_A$. Then
$$
\big(\D_K^mf\big)(A)=\lim_{s\to\be}\left(\sum_{j=0}^m(-1)^{m-j}\left(\begin{matrix}
m\\j\end{matrix}\right)f\big(A_s+jK\big)\right)
$$
for every $f\in\L_{\o,m}(\R)$, where the limit is taken in the strong operator topology.
\end{thm}

The proof of this theorem is similar to Theorem \ref{As}.

\begin{thm}
Let $m$ and $\o$ satisfy the hypotheses of Theorem {\rm\ref{oon}}
and let $f$ be a function in $\L_{\o,m}(\R)$ such that $|f(t)|\le\const(1+|t|)$.
Suppose that $A$ and $K$ are self-adjoint operators
with $\|K\|<\be$. Then $\cd_{f(A+jK)}\supset\cd_A$ for every $j$ and
$$
\big(\D_K^mf\big)(A)u=\sum_{j=0}^m(-1)^{m-j}\left(\begin{matrix}m\\j\end{matrix}\right)f\big(A+jK\big)u,\quad u\in\cd_A,
$$
where $\big(\D_K^mf\big)(A)$ is defined by {\rm\rf{defD4}}.
\end{thm}

The proof of this theorem is similar to the proof of Theorem \ref{1355}.

\

\

\noindent
\begin{tabular}{p{9cm}p{15cm}}
A.B. Aleksandrov & V.V. Peller \\
St-Petersburg Branch & Department of Mathematics \\
Steklov Institute of Mathematics  & Michigan State University \\
Fontanka 27, 191023 St-Petersburg & East Lansing, Michigan 48824\\
Russia&USA
\end{tabular}

\end{document}